\documentclass[twocolumn]{autart}    

\usepackage{graphicx}          
\usepackage{amsmath,amssymb,mathtools}
\usepackage{calrsfs}\DeclareMathAlphabet{\pazocal}{OMS}{zplm}{m}{n}
\usepackage{bm}
\usepackage{refcount}
\usepackage{color}
\usepackage{cite}
\usepackage[inline]{enumitem}
\usepackage{stackrel}
\usepackage{multirow}
\usepackage{arydshln}

\newcounter{example}
\newcounter{assumption}
\newcounter{theorem}
\newcounter{corollary}
\newcounter{definition}

\newcounter{lemma}

\newenvironment{smallarray}[1]
 {\null\,\vcenter\bgroup\tiny
  \arraycolsep=.13885em
  \hbox\bgroup$\array{@{}#1@{}}}
 {\endarray$\egroup\egroup\,\null}

\begin{document}	

\begin{frontmatter}

\title{Robust~Observer~Design~for Polytopic~Discrete-Time~Nonlinear~Descriptor~Systems\thanksref{footnoteinfo}} 

\thanks[footnoteinfo]{Corresponding author: T.~J.~Meijer.}

\author[TUe]{T.~J.~Meijer}\ead{t.j.meijer@tue.nl},    
\author[ASML]{V.~S.~Dolk}\ead{victor.dolk@asml.com},               
\author[TUe]{M.~S.~Chong}\ead{m.s.t.chong@tue.nl},  
\author[TUe]{W.~P.~M.~H.~Heemels}\ead{m.heemels@tue.nl}

\address[TUe]{Department of Mechanical Engineering, Eindhoven University of Technology, Eindhoven, The Netherlands}  
\address[ASML]{ASML, De Run 6665, 5504 DT Veldhoven, The Netherlands}             
          

\begin{abstract}                          
	This paper considers the design of robust state observers for a class of slope-restricted nonlinear descriptor systems with unknown time-varying parameters belonging to a known set. The proposed design accounts for process disturbances and measurement noise, while allowing for a trade-off between transient performance and sensitivity to noise and parameter mismatch. We exploit a polytopic structure of the system to derive linear-matrix-inequality-based synthesis conditions for robust parameter-dependent observers for the entire parameter set. In addition, we present (alternative) necessary and sufficient synthesis conditions for an important subclass within the considered class of systems and we show the effectiveness of the design for a numerical case study.
\end{abstract}
\end{frontmatter}

\section{Introduction}
In this paper, we propose a state observer for a class of slope-restricted nonlinear \textit{descriptor} systems with time-varying uncertain parameters that is robust in the sense that it is input-to-state stable (ISS)~\cite{Jiang2001} with respect to process disturbances, measurement noise and parameter mismatch. We consider a problem setting in which the true parameter values are unknown, however, they belong to a given compact set and a (possibly) time-varying estimate of the parameters (also in the set) is available. The model class, in particular with a parameter-dependent descriptor matrix, emerges when using implicit discretization schemes, even if the underlying continuous-time system does not feature a parameter-dependent  descriptor matrix.

The considered robust state estimation problem is relevant in its own right, for instance, to enable the synchronization of digital twins to their physical counterparts, see, e.g.,~\cite{Rasheed2020}. In addition, state observers are essential, for instance, in the context of output-based feedback stabilization, see, e.g.,~\cite{Heemels2010,Zemouche2016b}, and in certain sampling-based joint parameter and state estimation schemes, see, e.g.,~\cite{Meijer2021,Chong2015}. As such, the presented approach is relevant for a wide range of applications, especially, where the true parameters and/or the parameter estimates vary over time. Since the parameters belong to a known set, we aim to design observers that guarantee robustness for any possible (future) values of the time-varying true/estimated parameters. We will exploit a polytopic structure of the system, see also~\cite{Heemels2010,Zemouche2016b}, to develop LMI-based conditions to synthesize such observers off-line. The proposed observer is shown to guarantee ISS with respect to a general form of model mismatch. For some applications, e.g., joint parameter and state estimation~\cite{Meijer2021,Chong2015}, we are interested in ISS with respect to the parameter mismatch itself, which, as we will show, can also be guaranteed under (natural) additional assumptions. Observer designs have been proposed for descriptor systems, see, e.g.,~\cite{Estrada-Manzo2016,Hamdi2009}; for several classes of slope-restricted nonlinear (non-descriptor) systems, see, e.g.,~\cite{Zemouche2017,Ibrir2007,Chong2012}. These works consider nominal models, i.e., without uncertain parameters~\cite{Zemouche2017,Ibrir2007,Chong2012}, or assume that the exact parameters are known~\cite{Estrada-Manzo2016,Hamdi2009}. Our design considers three types of uncertainties, namely \begin{enumerate*}[label=(\alph*)] \item process disturbances, \item measurement noise and \item parameter mismatch, \end{enumerate*} and allows for a trade-off between transient performance and sensitivity to each type, which therefore differs from these existing works, see also Section~\ref{subsec:discussion}. Moreover, we propose alternative LMI-based synthesis conditions for an important special case of systems, with a parameter-independent descriptor matrix, which is also considered in, e.g.,~\cite{Bara2011,Chadli2008b}. For the latter class of systems, we show that the proposed conditions are not only sufficient, but also necessary thereby forming a strong indicator for the non-conservativeness of our design. To the best of our knowledge, no necessary and sufficient LMI-based synthesis conditions were available in the literature for this subclass of systems featuring uncertain parameters as well as slope-restricted nonlinearities. Interestingly, our results also reveal that for the considered subclass of systems, unlike other observer designs for polytopic systems in the literature, the structure obtained from our synthesis conditions does not lead to any additional conservatism compared to using observer gains which can depend arbitrarily on the parameters. 

The remainder of this paper is organized as follows. Section~\ref{sec:problem-definition} states the problem definition and Section~\ref{sec:motivation} motivates the considered model class. Section~\ref{sec:robust-observer-design} presents our observer design, both of the aforementioned synthesis conditions and a numerical example, which demonstrates the effectiveness of the design. Finally, Section~\ref{sec:conclusions} provides some conclusions. All proofs are provided in the Appendix.

\section{Preliminaries}
The following notation will be used throughout this paper. The sets of real and natural numbers are denoted by $\mathbb{R}=\left(-\infty,\infty\right)$ and $\mathbb{N}=\left\{0,1,2,\hdots\right\}$, respectively. Moreover, we denote $\mathbb{R}_{\geqslant 0} = \left[0,\infty\right)$, $\mathbb{N}_{\geqslant n} = \left\{n,n+1,n+2,\hdots\right\}$ and $\mathbb{N}_{[n,m]} = \{n,n+1,\hdots,m\}$, $n,m\in\mathbb{N}$. The set of symmetric (positive definite) matrices of size $n\times n$ is denoted by $\mathbb{S}^n$ ($\mathbb{S}^{n}_{\succ 0}$). The notation $\left(u,v\right)$ stands for $[u^\top\, v^\top]^\top$, where $u\in\mathbb{R}^m$ and $v\in\mathbb{R}^{n}$ with $n,m\in\mathbb{N}$. Let $\bm{e}_i\in\mathbb{R}^{n}$, $i\in\mathbb{N}_{[1,n]}$, be the $i$-th $n$-dimensional unit vector. For a vector $x\in\mathbb{R}^{n}$, $\|x\| = \sqrt{x^\top x}$ denotes its Euclidean norm. $\mathbb{B}^n_r\coloneqq\{x\in\mathbb{R}^n|\|x\|\leqslant r\}$ is the set of vectors with their norm bounded by $r$ and $\operatorname{Co}\{p_1,\hdots,p_N\}\subset\mathbb{R}^{m}$ denotes the convex hull of the points $\{p_i\}_{i\in\mathbb{N}_{[1,N]}}$ with $p_i\in\mathbb{R}^m$, $i\in\mathbb{N}_{[1,N]}$. For a matrix $A$, let $\operatorname{He}(A) = A + A^\top$, $A^\top$ denotes its transpose, $A^{-\top}$ the inverse of its transpose and $\sigma_{\operatorname{min}}(A)$ denotes its smallest singular value. The symbol $\star$ completes a symmetric matrix, e.g., ${\tiny\begin{bmatrix}\begin{smallmatrix} A & \star\\ B & C\end{smallmatrix}\end{bmatrix}}$ means ${\tiny\begin{bmatrix}\begin{smallmatrix} A & B^\top\\ B & C\end{smallmatrix}\end{bmatrix}}$. $M\prec 0$, $M\preccurlyeq 0$ and $M\succ 0$ mean, respectively, that $M\in\mathbb{S}^n$ is negative definite, negative semi-definite and positive definite. Let $\|A\|\coloneqq \sup_{x\neq 0} \|Ax\|/\|x\|$ denote the 2-norm of $A$, which equals $\sqrt{\lambda_{\mathrm{max}}(A^\top A)}$ with $\lambda_{\mathrm{max}}(A)$ the largest eigenvalue of $A$. The symbol $I$ is an identity matrix of appropriate dimensions and $\operatorname{diag}(D_1,D_2,\hdots,D_n)$ is a block diagonal matrix with $n$ diagonal blocks $D_i$, $i\in\mathbb{N}_{[1,n]}$. Finally, $y\leftarrow x$ means substituting $x$ for $y$.

A continuous function $\alpha\colon\mathbb{R}_{\geqslant 0}\rightarrow\mathbb{R}_{\geqslant 0}$ is a $\pazocal{K}$-function ($\alpha\in\pazocal{K}$), if it is strictly increasing and $\alpha\left(0\right)=0$. A continuous function $\beta\colon\mathbb{R}_{\geqslant 0}\times\mathbb{R}_{\geqslant 0}\rightarrow\mathbb{R}_{\geqslant 0}$ is a $\pazocal{KL}$-function ($\beta\in\pazocal{KL}$), if $\beta\left(\cdot,s\right)\in\pazocal{K}$ for each $s\in\mathbb{R}_{\geqslant 0}$, $\beta\left(r,\cdot\right)$ is non-increasing and $\beta\left(r,s\right)\rightarrow 0$ as $s\rightarrow\infty$ for each $r\in\mathbb{R}_{\geqslant 0}$. For a function $x\colon\mathbb{N}\rightarrow\mathbb{R}^{n}$, let $x_k$ denote the function evaluated at time $k\in\mathbb{N}$, i.e., $x_k\coloneqq x(k)$. We sometimes use $x_+$ to refer to $x_{k+1}=x(k+1)$. Finally, for a sequence $\{x_k\}_{k\in\mathbb{N}}$ we denote $\|x\|_{\infty,k} \coloneqq \max_{j\in\mathbb{N}_{[0,k]}} \|x_j\|$ and we denote the space of all bounded sequences taking values in $\mathbb{R}^{n}$ with $n\in\mathbb{N}$ by $\ell^{n}_\infty \coloneqq \{\{x_k\}_{k\in\mathbb{N}}\,|\,x_k\in\mathbb{R}^n,\,k\in\mathbb{N}\,\text{and}\,\sup_{k\in\mathbb{N}}\|x_k\| < \infty\}$.

\section{Problem definition}\label{sec:problem-definition}
Consider a discrete-time nonlinear descriptor system
\begin{align}
	E(p_{k+1})x_{k+1} &= A(p_k)x_k + G(p_k)\phi(Hx_k)+\nonumber\\
	&\hspace*{1.45cm} B(p_k)u_k +F(p_k)v_k,\label{eq:system}\\
	y_k &= Cx_k + Dw_k,\nonumber
\end{align}%
where $x_k\in\mathbb{R}^{n_x}$, $u_k\in\mathbb{R}^{n_u}$ and $y_k\in\mathbb{R}^{n_y}$ denote the state, input and output at discrete time $k\in\mathbb{N}$. The process noise at time $k\in\mathbb{N}$ is denoted by $v_k\in\mathbb{R}^{n_v}$ and the measurement noise by $w_k\in\mathbb{R}^{n_w}$. The unknown parameter $p_k\in\mathbb{P}$, $k\in\mathbb{N}$, belongs to a known compact set $\mathbb{P}\subset\mathbb{R}^{n_p}$ and we have $E,A\colon\mathbb{P}\rightarrow\mathbb{R}^{n_x\times n_x}$, $G\colon\mathbb{P}\rightarrow\mathbb{R}^{n_x\times n_\phi}$, $B\colon\mathbb{P}\rightarrow\mathbb{R}^{n_x\times n_u}$, $F\colon\mathbb{P}\rightarrow\mathbb{R}^{n_x\times n_v}$, $H\in\mathbb{R}^{n_\phi\times n_x}$, $C\in\mathbb{R}^{n_y\times n_x}$, $D\in\mathbb{R}^{n_y\times n_w}$ and $\phi\colon\mathbb{R}^{n_\phi}\rightarrow\mathbb{R}^{n_{\phi}}$. We assume that $E(p)$ is non-singular for all $p\in\mathbb{P}$ and, thus, the system~\eqref{eq:system} admits a unique solution defined for all $k\in\mathbb{N}$. In Section~\ref{sec:motivation}, we motivate the consideration of a parameter-dependent descriptor matrix $E$ and, in particular, the dependence on the next parameter $p_{k+1}$.
\setcounter{thm}{\theassumption}\stepcounter{assumption}
\begin{assum}\label{ass:parameter}
	The system~\eqref{eq:system} can be expressed in polytopic form, i.e., there exist functions $\xi_i\colon\mathbb{P}\rightarrow\mathbb{R}_{\geqslant 0}$, $i\in\pazocal{N}\coloneqq\mathbb{N}_{[1,N]}$, $N\in\mathbb{N}_{\geqslant 1}$, such that the mapping $\xi\coloneqq(\xi_1,\hdots,\xi_{N})$ satisfies $\xi(\mathbb{P})\subset\{\mu\in\mathbb{R}_{\geqslant 0}^N\,|\,\sum_{i\in\pazocal{N}}\mu_i=1\}$, and matrices $E_i,A_i\in\mathbb{R}^{n_x\times n_x}$, $G_i,H_i^\top\in\mathbb{R}^{n_x\times n_{\phi}}$, $B_i\in\mathbb{R}^{n_x\times n_u}$ and $F_i\in\mathbb{R}^{n_x\times n_v}$, $i\in\pazocal{N}$, such that, for all $p\in\mathbb{P}$,
\begin{equation}
		\begin{bmatrix}\begin{smallmatrix}
			E(p) & A(p) & B(p) & F(p) & G(p)
		\end{smallmatrix}\end{bmatrix} = \sum_{\mathclap{i\in\pazocal{N}}}\xi_i(p)\begin{bmatrix}\begin{smallmatrix}
			E_i & A_i & B_i & F_i & G_i
		\end{smallmatrix}\end{bmatrix}.
		\label{eq:poly-rep}
	\end{equation}
\end{assum}
\noindent The above holds, e.g., if~\eqref{eq:system} depends affinely on the parameters and $\mathbb{P}$ is (embedded in) a convex polytope, see, e.g.,~\cite{Heemels2010} for a class of uncertain LPV systems. 
\setcounter{thm}{\theassumption}\stepcounter{assumption}
\begin{assum}\label{ass:nonlinearity-condition}
	The nonlinearity $\phi$ in~\eqref{eq:system} is slope-restricted, i.e., there exists $\Lambda\in\mathbb{S}^{n_\phi}_{\succ 0}$ such that
	\begin{equation}
		\begin{bmatrix}\begin{smallmatrix}
			\Phi(x,y)\\
			y-x
		\end{smallmatrix}\end{bmatrix}^\top\begin{bmatrix}\begin{smallmatrix}
			2I & \star\\
			-\Lambda & 0
		\end{smallmatrix}\end{bmatrix}\begin{bmatrix}\begin{smallmatrix}
			\Phi(x,y)\\
			y-x
		\end{smallmatrix}\end{bmatrix} \leqslant 0,\text{ for all }x,y\in\mathbb{R}^{n_\phi},
		\label{eq:prop-slope-restr}
	\end{equation}
	where $\Phi(x,y) = \phi(y)-\phi(x)$, $x,y\in\mathbb{R}^{n_\phi}$.
\end{assum}
\noindent If, in addition, $\phi(0)=0$, the slope restriction~\eqref{eq:prop-slope-restr} implies that the nonlinearity is sector-bounded~\cite{Turner2020}.

Given the measured output $y_k$, input $u_k$ and an estimate $\hat{p}_k$ of the unknown parameter $p_k$, our objective is to design an observer for the system~\eqref{eq:system} that produces, at time $k\in\mathbb{N}$, a state estimate $\hat{x}_k$, which is robust in the sense that the error $e_k=\hat{x}_k-x_k$ is ISS with respect to the process noise $v_k$, measurement noise $w_k$ and the model (and parameter) mismatch ($\tilde{p}_k=\hat{p}_k-p_k$). For the precise definition of ISS, see Definition~\ref{def:ISS} in Section~\ref{sec:robust-observer-design}.

\section{A motivating example}\label{sec:motivation}
We motivate the consideration of models of the form~\eqref{eq:system} by showing how they arise when discretizing underlying continuous-time models of the form
\begin{equation}
	E_c(p(t))\frac{\mathrm{d}}{\mathrm{d}t}x(t) = A_c(p(t))x(t),\label{eq:ct-system}
\end{equation}
where $p(t)\in\mathbb{P}$, $t\in\mathbb{R}_{\geqslant 0}$, denotes the continuous-time unknown time-varying parameter vector and $E_c,A_c\colon\mathbb{P}\rightarrow\mathbb{R}^{n_x\times n_x}$. For notational compactness, we removed the nonlinearity, input and noise in~\eqref{eq:ct-system}, but they can be included in the developments below, e.g., using a semi-implicit approach~\cite{Ascher1997}. Let $t_k=kT_s$, $k\in\mathbb{N}$, denote the time at the $k$-th sampling instance with sampling period $T_{s}\in\mathbb{R}_{\geqslant 0}$. We discretize~\eqref{eq:ct-system} using Tustin's method, which combines two approximations of $E(p(\tau_k))x(\tau_k)$:
\begin{align}
	E_c(p(\tau_k))x(\tau_k) &\approx E_c(p_k)x_k+\theta E_c(p_k)\dot{x}(t_k),\\
	E_c(p(\tau_k))x(\tau_k) &\approx E_c(p_{k+1})x_{k+1}-\theta E_c(p_{k+1})\dot{x}(t_{k+1}),\nonumber
\end{align}
where $\tau_k=(t_k+t_{k+1})/2$, $\theta=T_s/2$ and $x_k\approx x(t_k)$, $k\in\mathbb{N}$. Equating both approximations and evaluating~\eqref{eq:ct-system} at both time instances yields
\begin{equation}
	\left(E_c(p_{k+1})-\theta A_c(p_{k+1})\right)x_{k+1} = \left(E_c(p_k)+\theta A_c(p_k)\right)x_k.\label{eq:discretized}
\end{equation}
This brief derivation shows that discrete-time descriptor models of the form~\eqref{eq:system} arise naturally, when discretizing a parameter-dependent continuous-time model, and, in particular, motivates the dependence of the $E$-matrix in~\eqref{eq:system} on $p_{k+1}$. We stress that models with a parameter-dependent $E$-matrix are prominent in discrete time in the sense that, even if $E_c$ is parameter-independent, an implicit discretization scheme yields a parameter-dependent $E$-matrix in discrete time, see~\eqref{eq:discretized}. The $E$-matrix of the discrete-time system is non-singular (as required in this paper) under the assumption that $\theta$ does not coincide, for any $k\in\mathbb{N}$, with any generalized eigenvalue(s) of the pair $(A(p_k),E(p_k))$, i.e., $\det(E(p_k)-\theta A(p_k))\neq 0$.

From this point on we consider discrete-time models of the form~\eqref{eq:system}, which includes the nonlinearity, input and noise. As mentioned before, these terms can be incorporated, e.g., by combining the above with forward Euler for the additional terms in a so-called semi-implicit approach, see, e.g.,~\cite{Ascher1997}. 
\setcounter{thm}{\theexample}\stepcounter{example}
\begin{exmp}\label{example}
	By applying semi-implicit discretization, as discussed above, to~\cite[Example 1]{Ibrir2007} (with added noise) we obtain a system of the form~\eqref{eq:system} with $p_k = (p_k^{\langle 1\rangle},p_k^{\langle 2\rangle}) \in \mathbb{P}=[\underline{p}^{\langle 1\rangle},\overline{p}^{\langle 1\rangle}]\times[\underline{p}^{\langle 2\rangle},\overline{p}^{\langle 2\rangle}]$ and, for $p\in\mathbb{P}$,
	\begin{align}
			E(p) &= \begin{bmatrix}\begin{smallmatrix}
				1+p^{\langle 2\rangle}/2 & p^{\langle 2\rangle}/2-2p^{\langle 1\rangle}\\
				p^{\langle 2\rangle} & 1+p^{\langle 2\rangle}
			\end{smallmatrix}\end{bmatrix},\, \,H^\top=\begin{bmatrix}\begin{smallmatrix} 1\\ 1\end{smallmatrix}\end{bmatrix}\nonumber,\\
			A(p) &= \begin{bmatrix}\begin{smallmatrix}
				1-p^{\langle 2\rangle}/2 & 2p^{\langle 1\rangle}-p^{\langle 2\rangle}/2\\
				-p^{\langle 2\rangle}  & 1-p^{\langle 2\rangle}
			\end{smallmatrix}\end{bmatrix},\,G(p) = \begin{bmatrix}\begin{smallmatrix}
				p^{\langle 2\rangle}\\
				2p^{\langle 2\rangle}
			\end{smallmatrix}\end{bmatrix},\label{eq:matrices}\\
			B(p)&=F(p)=p^{\langle 1\rangle}\begin{bmatrix}\begin{smallmatrix}
				1\\
				1\end{smallmatrix}\end{bmatrix},C = \begin{bmatrix}\begin{smallmatrix}1 & 0\end{smallmatrix}\end{bmatrix},\,D = 1,\nonumber
	\end{align}
	 and $\phi(z)=\sin(z)+z$, $z\in\mathbb{R}$. The parameters are the (constant) sampling time $p_k^{\langle 1\rangle}=T_{s}\in[\underline{p}^{\langle 1\rangle},\overline{p}^{\langle 1\rangle}]=[9.5,10.5]\cdot 10^{-3}$ and the time-varying Lipschitz constant of the system $p_k^{\langle 2\rangle}=\gamma(t_k) T_{s}/2\in[\underline{p}^{\langle 2\rangle},\overline{p}^{\langle 2\rangle}]=[0.0475,0.0525]$, $k\in\mathbb{N}$.
	
	The parameter set $\mathbb{P}$ can be taken as a polytope with $N=4$ vertices $\nu_1 = (\underline{p}^{\langle 1\rangle},\underline{p}^{\langle 2\rangle})$, $\nu_2 = (\overline{p}^{\langle 1\rangle},\overline{p}^{\langle 2\rangle})$, $\nu_3 = (\overline{p}^{\langle 1\rangle},\underline{p}^{\langle 2\rangle})$ and $\nu_4 = (\underline{p}^{\langle 1\rangle},\overline{p}^{\langle 2\rangle})$. Subdivide $\mathbb{P}$ into two simplices, i.e., triangles, denoted $\Delta_i \coloneqq \operatorname{Co}\{\nu_i,\nu_3,\nu_4\}$, $i\in\{1,2\}$. The functions $\xi_i$ in Assumption~\ref{ass:parameter}, $i\in\pazocal{N}$, are obtained by transforming to barycentric coordinates within each simplex. Let $\bar{\Delta}_1 \coloneqq \Delta_1$ and $\bar{\Delta}_2\coloneqq\Delta_2\setminus\operatorname{Co}\{\nu_3,\nu_4\}$, such that $\bar{\Delta}_1\cap\bar{\Delta}_2=\{0\}$, to obtain piecewise linear functions
	\begin{align}
		\begin{bmatrix}\begin{smallmatrix}
			\xi_3(p)\\
			\xi_4(p)
		\end{smallmatrix}\end{bmatrix} &\coloneqq T_i^{-1}(p-\nu_i),\quad \text{for }p\in\bar{\Delta}_i,\nonumber\\
		\xi_i(p) &\coloneqq	\begin{cases}
			1-\xi_3(p)-\xi_4(p),\quad &p\in\bar{\Delta}_i,\\
			0,\quad &p\notin\bar{\Delta}_i,
		\end{cases}
	\end{align}
	with $T_i = [\nu_3-\nu_i\quad \nu_4-\nu_i]$, $i\in\{1,2\}$. Since $\mathbb{P}=\cup_{i\in\{1,2\}}\bar{\Delta}_i$, the functions $\xi_i(p)$, $i\in\pazocal{N}$, are defined for all $p\in\mathbb{P}$. Evaluate the matrix-valued functions in~\eqref{eq:matrices} at each vertex, e.g., $E_i = E(\nu_i)$, to obtain $E_i,A_i,G_i,B_i,F_i$, $i\in\pazocal{N}$.
\end{exmp}

Since $E(p)$ is non-singular for all $p\in\mathbb{P}$, it seems convenient to pre-multiply the dynamics in~\eqref{eq:system} with $E^{-1}(p)$ to avoid dealing with the descriptor structure. This may destroy the polytopic structure, however, resulting in a more challenging synthesis problem, see, e.g.,~\cite{Estrada-Manzo2016}. Consider, for instance, the matrix inverse, for all $p\in[0,1]$,
\begin{equation}
	E^{-1}(p) = \begin{bmatrix}\begin{smallmatrix}
	1-p & p\\
	p & 1-p
	\end{smallmatrix}\end{bmatrix}^{-1} = \frac{1}{1-2p}\begin{bmatrix}\begin{smallmatrix}
		1-p & -p\\
	 	-p & 1-p
	\end{smallmatrix}\end{bmatrix},
\end{equation}
which is not polytopic despite $E(p)$ being polytopic. We will see that preserving the polytopic structure allows for the derivation of LMI-based synthesis conditions.

\section{Robust polytopic observer design}\label{sec:robust-observer-design}
We assume that estimates $\hat{p}_k\in\mathbb{P}$, $k\in\mathbb{N}$, of the unknown parameter vector $p_k$ are available at least one step into the future, i.e., $\hat{p}_k$ and $\hat{p}_{k+1}$ are known/estimated at time $k\in\mathbb{N}$. Based on these estimates, we aim to design a state observer of the form
\begin{align}
		&E(\hat{p}_{k+1})\hat{x}_{k+1} = A(\hat{p}_k)\hat{x}_k - L(\hat{p}_{k+1},\hat{p}_k)(C\hat{x}_k-y_k)+\nonumber\\
		&\quad B(\hat{p}_k)u_k +G(\hat{p}_k)\phi(H\hat{x}_k - K(\hat{p}_{k+1},\hat{p}_k)(C\hat{x}_k-y_k)),\label{eq:observer}
\end{align}
where $\hat{x}_k\in\mathbb{R}^{n_x}$ denotes the state estimate at time $k\in\mathbb{N}$ and $L\colon\mathbb{P}\times\mathbb{P}\rightarrow\mathbb{R}^{n_x\times n_y}$ and $K\colon\mathbb{P}\times\mathbb{P}\rightarrow\mathbb{R}^{n_\phi\times n_y}$ are to be designed. We will construct these observer gains in terms of polytopic functions, i.e., for all $\hat{p}_+,\hat{p}\in\mathbb{P}$,
\begin{align}
	L(\hat{p}_+,\hat{p}) &= \sum_{\mathclap{i,j\in\pazocal{N}}}\xi_i(\hat{p})\xi_j(\hat{p}_+)L_{ij},\label{eq:polytopic-L}\\
	K(\hat{p}_+,\hat{p}) &= \frac{Z(\hat{p}_+,\hat{p})}{\tau(\hat{p}_+,\hat{p})},\label{eq:K}
\end{align}
with
\begin{align}
	\tau(\hat{p}_+,\hat{p})=\sum_{\mathclap{i,j\in\pazocal{N}}}\xi_i(\hat{p})\xi_j(\hat{p}_+)\tau_{ij},\label{eq:tau}\\
	Z(\hat{p}_+,\hat{p})=\sum_{\mathclap{i,j\in\pazocal{N}}}\xi_i(\hat{p})\xi_j(\hat{p}_+)Z_{ij},\label{eq:Z}
\end{align}
where $\tau(\hat{p}_+,\hat{p})\in\mathbb{R}_{>0}$ for all $\hat{p}_+,\hat{p}\in\mathbb{P}$, such that they can be constructed by synthesizing a finite number of variables $L_{ij}\in\mathbb{R}^{n_x\times n_y}$, $\tau_{ij}\in\mathbb{R}_{>0}$ and $Z_{ij}\in\mathbb{R}^{n_\phi\times n_y}$, $i,j\in\pazocal{N}$. Here, $\xi_i$ are the functions in Assumption~\ref{ass:parameter}. The observer gain $L$ is polytopic in both $\hat{p}$ and $\hat{p}_+$ here. The observer gain $K$ in~\eqref{eq:K} with~\eqref{eq:tau}-\eqref{eq:Z} is not necessarily polytopic, however, by taking $\tau(\hat{p}_+,\hat{p})$ to be constant, i.e., $\tau(\hat{p}_+,\hat{p})=\bar{\tau}$ for some $\bar{\tau}\in\mathbb{R}_{>0}$ for all $\hat{p}_+,\hat{p}\in\mathbb{P}$, any polytopic observer gain can be constructed. Hence, this observer design is more general than when we restrict $K$ to be polytopic a priori as in, e.g.,~\cite{Heemels2010,Estrada-Manzo2016,Chadli2008b}. In Section~\ref{sec:nec-suff}, we will deviate from the structure in~\eqref{eq:observer} and~\eqref{eq:polytopic-L} by designing an observer with $L$, which is polytopic in $\hat{p}$, but not necessarily in $\hat{p}_+$ (the dependence may be of a general nature). This introduces additional flexibility in $L$ which we exploit to show that the proposed observer design conditions are necessary and sufficient for an important subclass of systems with a constant $E$.

Below we present conditions to synthesize observer gains as in~\eqref{eq:polytopic-L}-\eqref{eq:Z} that guarantee robust stability, where we take the following procedure: First, we consider robustness in terms of ISS with respect to process disturbances and measurement noise as well as a model mismatch induced by the parameter mismatch, which is often sufficient for output-based feedback control, see, e.g.,~\cite{Heemels2010}. In fact, we provide sufficient design conditions for the observer~\eqref{eq:observer}-\eqref{eq:Z} and general system setup~\eqref{eq:system}, in Section~\ref{sec:robust-mm}. Moreover, we present, in Section~\ref{sec:nec-suff}, an alternative set of observer design conditions that is not only sufficient, but also necessary for the case with a constant descriptor matrix thereby forming a strong indicator for the non-conservativeness of our design. Second, in Section~\ref{sub:param-mm}, we show that, under extra assumptions, an ISS property with respect to the parameter mismatch itself can be obtained for observers synthesized using either of the proposed conditions.

\subsection{Robustness with respect to model mismatch}\label{sec:robust-mm}
\sloppy Let $e_k\coloneqq \hat{x}_k-x_k$ be the state estimation error corresponding to~\eqref{eq:system} and~\eqref{eq:observer}. We group all terms that constitute the model mismatch (induced by $\tilde{p}_k\coloneqq \hat{p}_k-p_k$), and denote, for compactness, $\psi_k=\psi(\hat{p}_{k+1},\hat{p}_k,p_{k+1},p_k,x_{k+1},x_k,u_k,v_k)\coloneqq -(E(\hat{p}_{k+1})-E(p_{k+1}))x_{k+1} + (A(\hat{p}_k)-A(p_k))x_k + (B(\hat{p}_k)-B(p_k))u_k + (G(\hat{p}_k)-G(p_k))\phi(Hx_k)+(F(\hat{p}_k)-F(p_k))v_k$. We also denote $\tilde{\phi}_k=\tilde{\phi}(\hat{p}_{k+1},\hat{p}_k,e_k,x_k,w_k)\coloneqq \phi\big((H-K(\hat{p}_{k+1},\hat{p}_k)C)e_k+K(\hat{p}_{k+1},\hat{p}_k)Dw_k+Hx_k\big)-\phi(Hx_k)$. The error dynamics are then described by
\begin{align}
	&E(\hat{p}_{k+1})e_{k+1} = (A(\hat{p}_k)-L(\hat{p}_{k+1},\hat{p}_k)C)e_k + G(\hat{p}_k)\tilde{\phi}_k-\nonumber\\
	&\hspace*{2.4cm} F(\hat{p}_k)v_k +L(\hat{p}_{k+1},\hat{p}_k)Dw_k+ \psi_k.
	\label{eq:error-system}
\end{align}
Next, we give the definition of ISS~\cite{Jiang2001}, which incorporates our aim to achieve robustness for the entire parameter set.
\setcounter{thm}{\thedefinition}\stepcounter{definition}
\begin{defn}\label{def:ISS}
	The system~\eqref{eq:error-system} is said to be ISS with respect to $v$, $w$ and $\psi$, if there exist functions $\beta\in\pazocal{KL}$ and $\gamma_v,\gamma_w,\gamma_{\psi}\in\pazocal{K}$ such that, for all $e_0\in\mathbb{R}^{n_x}$, $v\in\ell_\infty^{n_v}$, $w\in\ell_{\infty}^{n_w}$, $\psi\in\ell_{\infty}^{n_x}$, $\{p_k\}_{k\in\mathbb{N}}$ and $\{\hat{p}_k\}_{k\in\mathbb{N}}$ with $p_k,\hat{p}_k\in\mathbb{P}$,
	\begin{align}
		\|e_k\|&\leqslant \beta(\|e_0\|,k)+\gamma_v(\|v\|_{\infty,k-1})+		\label{eq:iss}\\
		&\qquad\gamma_w(\|w\|_{\infty,k-1})+\gamma_\psi(\|\psi\|_{\infty,k-1}),\quad k\in\mathbb{N}.\nonumber
	\end{align}	
\end{defn}
\noindent To synthesize observers~\eqref{eq:observer} that render the error system~\eqref{eq:error-system} ISS with respect to disturbances, measurement noise and model mismatch, we construct suitable ISS-Lyapunov functions~\cite{Jiang2001}. Inspired by~\cite{Heemels2010,Daafouz2001}, we focus on a class of polytopic parameter-dependent quadratic ISS-Lyapunov functions, which leads to a systematic design procedure.
\setcounter{thm}{\thedefinition}\stepcounter{definition}
\begin{defn}\label{def:poly-ISS}
	The system~\eqref{eq:error-system} is said to be absolutely poly-quadratically ISS with respect to $v$, $w$ and $\psi$, if it admits a poly-quadratic ISS-Lyapunov function $V\colon\mathbb{R}^{n_x}\times\mathbb{P}\rightarrow\mathbb{R}_{\geqslant 0}$ of the form $V(e,\hat{p})=\sum_{i\in\pazocal{N}}e^\top\xi_i(\hat{p})P_ie$ (with $P_i\in\mathbb{S}^{n_x}$, $i\in\pazocal{N}$), satisfying\footnote{Without loss of generality, we can obtain~\eqref{eq:descent-cond} by scaling $V$ to get a coefficient $-1$ in $-\|e\|^2$. The lower bound in~\eqref{eq:lyap-bounds} is guaranteed by~\eqref{eq:descent-cond}, since $V(e_+,\hat{p}_+)\geqslant 0$ and~\eqref{eq:descent-cond} holds with $\|v\|=\|w\|=\|\psi\|=0$, and, hence, $a\geqslant 1$.}
	\begin{align}
		\|e\|^2 \leqslant V(e,\hat{p}) &\leqslant a\|e\|^2,\label{eq:lyap-bounds}\\
		V(e_+,\hat{p}_+) - V(e,\hat{p}) &\leqslant -\|e\|^2 + \kappa_v\|v\|^2 + \nonumber\\
		&\qquad \kappa_w\|w\|^2 + \kappa_{\psi}\|\psi\|^2,\label{eq:descent-cond}
	\end{align}	
	for all $\tilde{\phi}(\hat{p}_+,\hat{p},e,x,w)=\phi((H-K(\hat{p}_+,\hat{p})C)e+K(\hat{p}_+,\hat{p})Dw+Hx)-\phi(Hx)$ with $\phi$ satisfying Assumption~\ref{ass:nonlinearity-condition}, $\hat{p}_+,\hat{p},p_+,p\in\mathbb{P}$, $e,x_+,x,\psi\in\mathbb{R}^{n_x}$, $u\in\mathbb{R}^{n_u}$, $v\in\mathbb{R}^{n_v}$ and $w\in\mathbb{R}^{n_w}$, with $a\in\mathbb{R}_{\geqslant 1}$, $\kappa_v,\kappa_w,\kappa_{\psi}\in\mathbb{R}_{>0}$ and $e_+$ satisfying $E(\hat{p}_+)e_+=(A(\hat{p})-L(\hat{p}_+,\hat{p})C)e+G(\hat{p})\tilde{\phi}(\hat{p}_+,\hat{p},e,x,w)-F(\hat{p})v+L(\hat{p}_+,\hat{p})Dw+\psi(\hat{p}_+,\hat{p},p_+,p,x_+,x,u,v)$.
\end{defn}
\noindent The existence of a poly-quadratic ISS-Lyapunov function for~\eqref{eq:error-system} implies that the error system~\eqref{eq:error-system} is ISS, see~\cite{Heemels2010}. By constructing poly-quadratic ISS Lyapunov functions, we derive sufficient LMI-based conditions for~\eqref{eq:error-system} to be (absolutely) poly-quadratically ISS, as stated in the theorem below.
\setcounter{thm}{\thetheorem}\stepcounter{theorem}
\begin{thm}\label{thm:observer-lmi}
	Consider system~\eqref{eq:system} and observer~\eqref{eq:observer} satisfying Assumptions~\ref{ass:parameter} and~\ref{ass:nonlinearity-condition}. Let $\underline{\sigma} \coloneqq \min_{p\in\mathbb{P}}\sigma_{\mathrm{min}}(E(p))$ and suppose there exist matrices $P_i\in\mathbb{S}^{n_x}$, $X_{i}\in\mathbb{R}^{n_x\times n_x}$, $Y_{ij}\in\mathbb{R}^{n_x\times n_y}$, $Z_{ij}\in\mathbb{R}^{n_\phi\times n_y}$ and scalars $\tau_{ij},\kappa_v,\kappa_w,\kappa_{\psi}\in\mathbb{R}_{>0}$, $i,j\in\pazocal{N}$, such that
	\begin{align}
		\underbrace{\begin{bmatrix}\begin{smallmatrix}
			\operatorname{He}(X_{i}E_j)-P_j & \star & \star & \star & \star & \star\\
			(X_{i} A_i-Y_{ij}C)^\top & P_i-I & \star & \star & \star & \star\\
			-(X_{i}G_i)^\top & \Lambda(\tau_{ij} H-Z_{ij}C) & 2\tau_{ij} I & \star & \star & \star\\
			-(X_{i}F_i)^\top & 0 & 0 & \kappa_vI & \star & \star\\
			(Y_{ij}D)^\top & 0 & (\Lambda Z_{ij}D)^\top & 0 & \kappa_wI & \star\\
			X_{i}^\top & 0 & 0 & 0 & 0 & \kappa_\psi I
		\end{smallmatrix}\end{bmatrix}}_{\eqqcolon M_{ij}}{\tiny\succ 0,}
		\label{eq:lmi}
	\end{align}
	for all $i,j\in\pazocal{N}$. Then, the matrices $X_i$ are non-singular, $P_i\succ 0$, $i\in\pazocal{N}$, and the state estimation error system~\eqref{eq:error-system} for system~\eqref{eq:system} and observer~\eqref{eq:observer} with~\eqref{eq:polytopic-L}-\eqref{eq:Z} and $L_{ij}=X^{-1}_{i}Y_{ij}$, $i,j\in\pazocal{N}$, is (absolutely) poly-quadratically ISS with respect to $v$, $w$ and $\psi$ with $\beta(s,r)\coloneqq\sqrt{\kappa_{\psi}}\underline{\sigma}\sqrt{\rho^r}s$, for some $\rho\in(0,1)$, $\gamma_v(s)\coloneqq\sqrt{\kappa_{\psi}\kappa_v}\underline{\sigma}s$, $\gamma_w(s)\coloneqq\sqrt{\kappa_{\psi}\kappa_w}\underline{\sigma}s$ and $\gamma_\psi(s)\coloneqq\kappa_\psi\underline{\sigma}s$, $r,s\in\mathbb{R}_{\geqslant 0}$. Moreover, $V(e,\hat{p})=\sum_{i\in\pazocal{N}}e^\top\xi_i(\hat{p})P_ie$ is a poly-quadratic ISS-Lyapunov function for~\eqref{eq:error-system} satisfying~\eqref{eq:lyap-bounds}-\eqref{eq:descent-cond} with $a=\kappa_\psi\underline{\sigma}^{2}>1$.
\end{thm}
\noindent Theorem~\ref{thm:observer-lmi} provides an LMI condition, with dimensions that scale with $N^2$ (recall that $N$ is the number of vertices of the polytope from Assumption~\ref{ass:parameter}), which enables us to synthesize observer gains $L(\hat{p}_+,\hat{p})$ and $K(\hat{p}_+,\hat{p})$ such that~\eqref{eq:error-system} is ISS. When the parameter estimate is constant, i.e., $\hat{p}_{k+1}=\hat{p}_k$, the condition can be relaxed by setting $j=i$ in~\eqref{eq:lmi}, such that the dimensions of the LMI scale linearly in $N$. The ISS-gains can, for example, be minimized via $\kappa_v$, $\kappa_w$ and $\kappa_\psi$ by the optimization regime
\begin{equation}
	\begin{aligned}
		&\text{minimize}  & & c_v\kappa_v+c_w\kappa_w+c_\psi\kappa_{\psi},\\
		&\text{subject to} & & \eqref{eq:lmi},
	\end{aligned}
\end{equation}
where $c_v,c_w,c_\psi\in\mathbb{R}_{\geqslant 0}$ with $c_v+c_w+c_{\psi}=1$. By tuning the weights $c_v,c_w,c_\psi$, the proposed conditions allow us to make trade-offs between the transient performance, sensitivity to the different noise levels and model mismatch. In the nominal case ($\tilde{p}_k=0$, $v_k=0$, $w_k=0$, $k\in\mathbb{N}$), the observer~\eqref{eq:observer} with the obtained polytopic gains renders the error dynamics~\eqref{eq:error-system} globally exponentially stable. 

\subsubsection{Discussion of Theorem~\ref{thm:observer-lmi} relative to other works}\label{subsec:discussion}
We discuss Theorem~\ref{thm:observer-lmi} in relation to relevant existing results regarding two aspects: \begin{enumerate*}[label=(\alph*)] \item the class of nonlinearities $\phi$ and \item the conservatism of the LMI condition~\eqref{eq:lmi}\end{enumerate*}. Starting with (a): For systems with a scalar nonlinearity, i.e., where $n_\phi=1$, the variables $\tau_{ij}$ fulfill the same role as the diagonal multiplier matrix introduced in~\cite{Chong2012}. The class of vector-valued ($n_\phi>1$) slope--restricted nonlinearities considered in Assumption~\ref{ass:nonlinearity-condition} is very general and guaranteeing absolute stability with respect to such a general class of nonlinearities may for many applications yield conservative results. As a result, Theorem~\ref{thm:observer-lmi} applies to a wide range of nonlinear systems, however, it typically pays off to exploit specific properties/structure of the nonlinearity if such information is available. For instance, if the entries of the vector-valued nonlinearity are decoupled, i.e., $\phi_i(z)$ can be written as $\phi_i(z_i)$ for $\phi_i\colon\mathbb{R}\rightarrow\mathbb{R}$, $i\in\mathbb{N}_{[1,n_\phi]}$, we can deal with each decoupled nonlinearity separately by associating each with its own scalar $\tau_{k,ij}\in\mathbb{R}_{>0}$, $k\in\mathbb{N}_{[1,n_\phi]}$, $i,j\in\pazocal{N}$, which leads to a diagonal multiplier matrix $\pazocal{T}_{ij}=\operatorname{diag}\{\tau_{1,ij},\hdots,\tau_{n_{\phi},ij}\}$ such as in~\cite{Chong2012}. Other relaxations for specific classes of nonlinearities that may reduce conservatism are proposed in, for instance,~\cite{Zemouche2017}.

Regarding (b): In the noiseless linear non-descriptor case ($\phi=0$, $v=0$, $w=0$ and $E(p)=I$ for all $p\in\mathbb{P}$),~\eqref{eq:lmi} reduces to the conditions in~\cite[Theorem 2]{Heemels2010}, which were shown to also be necessary in~\cite{Pandey2018}. In the next section, we restrict the considered class of systems to feature a constant $E$-matrix, which allows us to introduce more flexibility in how the observer gains depend on the parameters as well as introduce more slack variables. As a result, we can provide necessary and sufficient LMI-based synthesis conditions in this setting. The necessity indicates the non-conservativeness of our synthesis conditions. The conditions in Theorem~\ref{thm:observer-lmi} can be relaxed by applying the dilated LMI techniques used in, e.g.,~\cite{Oliveira2001,Ebihara2005,Bara2011}, which would allow us to introduce slack variables $X_{ij}$, $i,j\in\pazocal{N}$, (instead of just $X_i$ as in~\eqref{eq:lmi}) thereby reducing conservativeness. 

\subsection{Necessary and sufficient conditions for the case $E(p)=\bar{E}$ for all $p\in\mathbb{P}$}\label{sec:nec-suff}
Next, we provide necessary and sufficient conditions to synthesize observers of the form~\eqref{eq:observer}, that guarantee absolute poly-quadratic ISS, for systems~\eqref{eq:system} with a constant descriptor matrix, as formalized below.
\setcounter{thm}{\theassumption}\stepcounter{assumption}
\begin{assum}\label{ass:Ebar}
	It holds that $E(p)=\bar{E}$, for all $p\in\mathbb{P}$, for some non-singular $\bar{E}\in\mathbb{R}^{n_x\times n_x}$.
\end{assum}
\noindent A significant class of systems satisfies Assumption~\ref{ass:Ebar} including, for instance, standard state-space models, but also models obtained through forward Euler discretization of some underlying continuous-time descriptor system with a parameter-independent $E$-matrix. Note that $\bar{E}$ is still non-singular and, hence, pre-multiplying the dynamics by $\bar{E}^{-1}$ yields a non-descriptor system. However, compared to existing works, e.g.,~\cite{Zemouche2016b,Pandey2018,Heemels2010}, we still show necessity for a more general class of nonlinear systems for which, to the best of our knowledge, such results were not available in the literature. In addition, our results reveal that the obtained observer gains are no more conservative than observer gains which arbitrarily depend on both $\hat{p}_{k+1}$ as well as $\hat{p}_k$, as detailed below.

For the above class of systems, we synthesize observers~\eqref{eq:observer} with arbitrary observer gains $L\colon\mathbb{P}\times\mathbb{P}\rightarrow\mathbb{R}^{n_x\times n_y}$ and $K\colon\mathbb{P}\times\mathbb{P}\rightarrow\mathbb{R}^{n_\phi\times n_y}$. In other words, we do not impose any restrictions on how $L$ and $K$ depend on $\hat{p}_k$ and $\hat{p}_{k+1}$. In doing so, we not only provide necessary and sufficient LMI-based synthesis conditions for the considered class of observers~\eqref{eq:observer}, but we will also uncover a necessary structure for the observer gains $L$ and $K$. Finally, we adopt the additional assumption on the functions $\xi_i$, $i\in\pazocal{N}$, from Assumption~\ref{ass:parameter} that $\{\bm{e}_i\}_{i\in\pazocal{N}}\subset\xi(\mathbb{P})$ (which holds for the functions $\xi_i$ in Example~\ref{example}). This assumption is mild/natural when studying necessary stability conditions for uncertain systems, since it means that we are not overapproximating the parameter set in order to satisfy Assumption~\ref{ass:parameter} which inherently leads to conservatism. In the case where we do have to over approximate the parameter set and, thereby, introduce conservatism, we can redefine $\mathbb{P}$ to satisfy this assumption and apply the results below to guarantee necessity for overapproximating the parameter set.
\setcounter{thm}{\thetheorem}\stepcounter{theorem}
\begin{thm}\label{thm:suffnec}
	Consider system~\eqref{eq:system} satisfying Assumptions~\ref{ass:parameter}-\ref{ass:Ebar} with $\{\bm{e}_i\}_{i\in\pazocal{N}}\subset\xi(\mathbb{P})$. There exists an observer of the form~\eqref{eq:observer}, with some observer gains $L\colon\mathbb{P}\times\mathbb{P}\rightarrow\mathbb{R}^{n_x\times n_y}$, $K\colon\mathbb{P}\times\mathbb{P}\rightarrow\mathbb{R}^{n_\phi\times n_y}$, that renders the error system~\eqref{eq:error-system} absolutely poly-quadratically ISS with respect to $v$, $w$ and $\psi$, if and only if there exist matrices $P_i\in\mathbb{S}^{n_x},X_{ij}\in\mathbb{R}^{n_x\times n_x},Y_{ij}\in\mathbb{R}^{n_x\times n_y}$, $Z_{ij}\in\mathbb{R}^{n_\phi\times n_y}$ and scalars $\tau_{ij},\kappa_v,\kappa_w,\kappa_\psi\in\mathbb{R}_{>0}$, $i,j\in\pazocal{N}$, such that
	\begin{equation}
		\underbrace{\begin{bmatrix}\begin{smallmatrix}
			\operatorname{He}(X_{ij}\bar{E})-P_j & \star & \star & \star & \star & \star\\
			(X_{ij}A_i-Y_{ij}C)^\top & P_i-I & \star & \star & \star & \star\\
			-(X_{ij}G_i)^\top & \Lambda(\tau_{ij}H - Z_{ij}C) & 2\tau_{ij}I & \star & \star & \star\\
			-(X_{ij}F_i)^\top & 0 & 0 & \kappa_vI & \star & \star\\
			(Y_{ij}D)^\top & 0 & (\Lambda Z_{ij}D)^\top & 0 & \kappa_wI & \star\\
			X_{ij}^\top & 0 & 0 & 0 & 0 & \kappa_{\psi}I
		\end{smallmatrix}\end{bmatrix}}_{\eqqcolon N_{ij}} \succ 0,
		\label{eq:nec-cond}
	\end{equation}
	for all $i,j\in\pazocal{N}$. In fact, then, $P_i\succ 0$, the matrices $X_{ij}$ are non-singular, $i,j\in\pazocal{N}$, and the observer~\eqref{eq:observer} with 
	\begin{equation}
		L(\hat{p}_+,\hat{p})=\sum_{\mathclap{i\in\pazocal{N}}}\xi_i(\hat{p})X_i^{-1}(\hat{p}_+)Y_i(\hat{p}_+),\,\hat{p}_+,\hat{p}\in\mathbb{P},	
		\label{eq:L}
	\end{equation}
	where
	\begin{equation}
		\begin{bmatrix}\begin{smallmatrix} X_i(\hat{p}_+) & Y_i(\hat{p}_+)\end{smallmatrix}\end{bmatrix} = \sum_{\mathclap{j\in\pazocal{N}}}\xi_j(\hat{p}_+)\begin{bmatrix}\begin{smallmatrix} X_{ij} & Y_{ij}\end{smallmatrix}\end{bmatrix},
		\label{eq:XY}
	\end{equation}
	and $K$ as defined in~\eqref{eq:K}-\eqref{eq:Z} renders the error system~\eqref{eq:error-system} absolutely poly-quadratically ISS with respect to $v$, $w$ and $\psi$ with $\beta$, $\gamma_v$, $\gamma_w$ and $\gamma_\psi$ as defined in Theorem~\ref{thm:observer-lmi}. Moreover, $V(e,\hat{p})=\sum_{i\in\pazocal{N}}e^\top \xi_i(\hat{p})P_ie$ is a poly-quadratic ISS-Lyapunov function for~\eqref{eq:error-system} satisfying~\eqref{eq:lyap-bounds}-\eqref{eq:descent-cond} with $a=\kappa_{\psi}\underline{\sigma}^2>1$ and $\underline{\sigma}\coloneqq \sigma_{\mathrm{min}}(\bar{E})$.
\end{thm}
\noindent The conditions in Theorem~\ref{thm:suffnec} are LMIs in the decision variables $X_{ij}$, $Y_{ij}$, $P_i$, $Z_{ij}$, $\tau_{ij}$, $\kappa_v$, $\kappa_w$ and $\kappa_{\psi}$. Since they are necessary and sufficient, these conditions can be used to synthesize observers of the form~\eqref{eq:observer} without introducing additional conservatism. In fact, Theorem~\ref{thm:suffnec} shows that restricting $L$ and $K$ to be of the form in~\eqref{eq:L}-\eqref{eq:XY} and~\eqref{eq:K}-\eqref{eq:Z}, respectively, also does not introduce any further conservatism. 

As before, the synthesized observer gain $K$ satisfies~\eqref{eq:K}-\eqref{eq:Z}, however, considering systems~\eqref{eq:system} with a constant descriptor matrix allowed us to derive LMIs for synthesizing observers with a more general $L$ as given in~\eqref{eq:L}-\eqref{eq:XY}, which is polytopic in $\hat{p}$, but no longer in $\hat{p}_+$. To see that~\eqref{eq:L} with~\eqref{eq:XY} is indeed more general than~\eqref{eq:polytopic-L}, note that by setting $X_{ij}=X_i$, for all $j\in\pazocal{N}$, the structure from Theorem~\ref{thm:observer-lmi} is recovered. Note that, in general, it is not possible to invert $X_i(\hat{p}_+)$ analytically and, hence, this inverse needs to be computed on-line. In some applications it may be desirable to set $X_{ij}=X_i$, for all $j\in\pazocal{N}$, thereby, avoiding the need to invert $X_i(\hat{p}_+)$ on-line and, hence, reducing the on-line computational burden at the cost of introducing conservatism.

\subsection{Robustness with respect to parameter mismatch}
\label{sub:param-mm}
Under some extra assumptions, ISS with respect to the model mismatch implies ISS with respect to the parameter mismatch, as formalized next.
\setcounter{thm}{\thecorollary}\stepcounter{corollary}
\begin{cor}\label{cor:param-mismatch-iss}
	Suppose that the conditions in Theorem~\ref{thm:observer-lmi} or Theorem~\ref{thm:suffnec} hold and that the functions $\xi_i$, $i\in\pazocal{N}$, are continuous. Given $K_x,K_u,K_v\in\mathbb{R}_{\geqslant 0}$, there exists a $\pazocal{K}$-function $\gamma_{\tilde{p}}$ such that, for all solutions $(x,u,v)$ to~\eqref{eq:system} that satisfy $x_k \in\mathbb{B}_{K_x}^{n_x}$, $u_k\in\mathbb{B}_{K_u}^{n_u}$ and $v_k\in\mathbb{B}_{K_v}^{n_v}$ and, for all $k\in\mathbb{N}$, 
	\begin{align}
		\|e_k\| &\leqslant \beta(\|e_0\|,k)+\gamma_v(\|v\|_{\infty,k-1})+\nonumber\\
		&\hfill \gamma_w(\|w\|_{\infty,k-1})+\gamma_{\tilde{p}}(\|\tilde{p}\|_{\infty,k-1}),
		\label{eq:iss-param}
	\end{align}
	for all $e_0\in\mathbb{R}^{n_x}$, $w\in\ell_{\infty}^{n_w}$, $\{p_k\}_{k\in\mathbb{N}}$ and $\{\hat{p}_k\}_{k\in\mathbb{N}}$ with $p_k,\hat{p}_k\in\mathbb{P}$, $k\in\mathbb{N}$. Here, $\beta,\gamma_v,\gamma_w$ are as in Theorem~\ref{thm:observer-lmi}.
\end{cor}
Corollary~\ref{cor:param-mismatch-iss} states that, if $u$ is bounded and keeps $x$ bounded for all admissible initial conditions $x_0$ and disturbances $v$, the error system~\eqref{eq:error-system} corresponding to system~\eqref{eq:system} with observer~\eqref{eq:observer} obtained from Theorem~\ref{thm:observer-lmi} or Theorem~\ref{thm:suffnec}, satisfies an ISS-like property with respect to $v$, $w$ and $\tilde{p}$.

\subsection{Illustrative example}
	\begin{figure}[!t]
		\centering
		\vspace*{-0.38cm}
		\includegraphics[width=.5\textwidth]{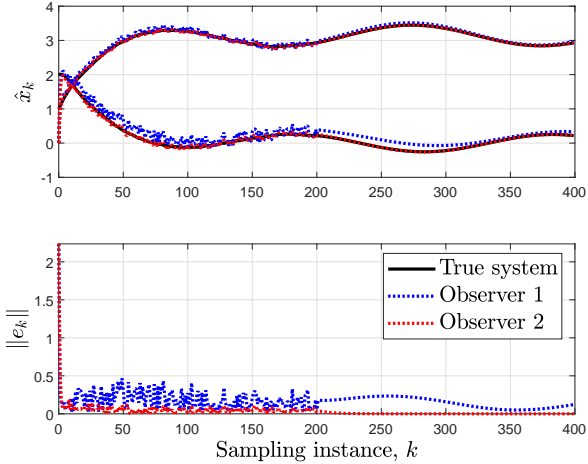}
		\caption{State estimates and corresponding estimation errors.}%
		\label{fig:results}%
	\end{figure}
	We demonstrate the effectiveness of our results using a case study based on Example~\ref{example}. Consider again the system~\eqref{eq:system} with~\eqref{eq:matrices} and its polytopic representation derived in Example~\ref{example}. The nonlinearity $\phi$ satisfies Assumption~\ref{ass:nonlinearity-condition} with $\Lambda=2$. Hence, the conditions of Theorem~\ref{thm:observer-lmi} are satisfied. The functions $\xi_i$, $i\in\{3,4\}$, are continuous at the boundary $\operatorname{Co}\{\nu_3,\nu_4\}$ and, thus, $\xi_i$, $i\in\pazocal{N}$, are continuous on $\mathbb{P}$. We minimize $(\kappa_v+5\kappa_w+0.01\kappa_\psi)/6.01$ subject to~\eqref{eq:lmi}. The optimal solution achieves $\kappa_v = 0.5258$, $\kappa_w = 6.6402$ and $\kappa_{\psi}=1.3535\cdot 10^3$. Computations show that $\underline{\sigma} = 1.01$. We implement two observers~\eqref{eq:observer}: One with the constant (inaccurate) parameter estimate $\hat{p}_k = \nu_1=(9.5\cdot 10^{-3},4.75\cdot 10^{-2})$ and the other with the exact parameter, i.e., $\hat{p}_k=p_k$, $k\in\mathbb{N}$. We apply both observers to the system~\eqref{eq:system} with $p_k = (0.01,0.11+0.11\sin(0.001\pi k))$, $u_k = 3\sin(0.01\pi k)$. On the time interval $k\in\mathbb{N}_{[0,200]}$, $v_k$ is a uniformly distributed random variable in the interval $[-1,1]$ and the entries of $w_k$ are uniformly distributed random variables in the interval $[-0.1,0.1]$. For $k\in\mathbb{N}_{\geqslant 201}$, no more noise is present, i.e., $v_k=0$ and $w_k=0$ for $k\in\mathbb{N}_{\geqslant 201}$. We simulate system~\eqref{eq:system} with $x_0 = (1,2)$ and both observers~\eqref{eq:observer} with $\hat{x}_0 = (0,0)$. 
	
	Fig.~\ref{fig:results} shows that the state estimates initially converge to a neighbourhood of the true states. We also see that after $k=200$ (when the noise is absent), the state estimate of observer 2 converges asymptotically to the true state whereas observer 1 retains some error due to the parameter mismatch. This behaviour nicely illustrates that the state estimation error system is ISS with respect to the process disturbances, measurement noise and model/parameter mismatch. We can verify, using Lyapunov-based tools~\cite{Khalil1996} and the results in Corollary~\ref{cor:param-mismatch-iss}, that ISS with respect to parameter mismatch is achieved.

\section{Conclusions}\label{sec:conclusions}
We derived LMI-based conditions to synthesize polytopic observers for a class of \textit{nonlinear descriptor} systems with uncertain parameters. The proposed approach leads to observers that can be synthesized a priori and which guarantee ISS with respect to process disturbances, measurement noise and model/parameter mismatch for the entire parameter set. This is particularly useful in applications where the parameter estimates are time-varying, such as in the joint parameter and state estimation scheme using a multi-observer in~\cite{Meijer2021,Chong2015}. The presented LMI conditions allow trade-offs between transient performance and sensitivity to noise and model/parameter mismatch. For the class of systems with parameter-independent $E$-matrix, we show that we can exploit additional flexibility in the structure of the observer gains to derive alternative conditions that are both necessary and sufficient, being an indicator for the non-conservativeness of our synthesis conditions. For the latter class of systems, we also show that the obtained structure for the observer gains does not lead to additional conservatism with respect to using observer gains which depend arbitrarily on the parameters. Finally, the proposed conditions were used to synthesize observers for a numerical example, showing their strengths.

\section*{Appendix}
First, we introduce some preliminaries which we use in proving our results.
	\setcounter{thm}{\thelemma}\stepcounter{lemma}
\begin{lem}[{Finsler's lemma~\cite[Lemma 3]{Ishihara2017}}]\label{lem:finsler}
	Let $\mathbb{P}\subseteq\mathbb{R}^{d}$, $Q\colon\mathbb{P}\rightarrow\mathbb{S}^n$ and $B\colon\mathbb{P}\rightarrow\mathbb{R}^{m\times n}$. Then, the following statements are equivalent:
	\begin{enumerate}[label=\alph*)]
		\item For all $p\in\mathbb{P}$, it holds that $x^\top Q(p)x<0$ for all $x\in\mathbb{R}^{n_x}\setminus\{0\}$ such that $B(p)x=0$.
		\item For all $p\in\mathbb{P}$, there exists $\mu(p)\in\mathbb{R}$ such that $Q(p)-\mu(p)B^\top(p)B(p)\prec 0$.
	\end{enumerate}
\end{lem}
	\setcounter{thm}{\thelemma}\stepcounter{lemma}
\begin{lem}[{Projection lemma~\cite[Lemma 4.15]{Scherer2011}}]\label{lem:proj}
	Let $Q\in\mathbb{S}^n$, $A\in\mathbb{R}^{m\times n}$ and $B\in\mathbb{R}^{p\times n}$. There exists $X\in\mathbb{R}^{m\times p}$ such that
	\begin{equation}
		Q+A^\top XB+B^\top X^\top A\succ 0,
	\end{equation}
	if and only if
	\begin{equation}
		A_\perp^\top QA_{\perp}\succ 0,\text{ and }B_{\perp}^\top QB_{\perp}\succ 0,
	\end{equation}
	where $A_{\perp}$ and $B_{\perp}$ denote arbitrary matrices whose columns form a basis of $\ker A$ and $\ker B$, respectively.
\end{lem}
	\setcounter{thm}{\thelemma}\stepcounter{lemma}
\begin{lem}[(Lossless) $\pazocal{S}$-procedure~\cite{Boyd1994}]\label{lem:S-proc}
	Let $Q_i\in\mathbb{S}^{n}$, $i\in\mathbb{N}_{[0,N]}$. Suppose that
	\begin{equation}
		\exists \tau_i\in\mathbb{R}_{\geqslant 0},\,i\in\mathbb{N}_{[1,N]},\text{ such that }Q_0-\sum_{\mathclap{i\in\mathbb{N}_{[1,N]}}}\tau_iQ_i\succ 0,\label{eq:S1}
	\end{equation}
	then, it holds that
	\begin{align}
		&x^\top Q_0x>0,\text{ for all }x\in\mathbb{R}^{n}\setminus\{0\}\text{ such that }\nonumber\\
		&\qquad x^\top Q_ix\geqslant 0\text{ for all }i\in\mathbb{N}_{[1,N]}.\label{eq:S2}
	\end{align}
	Conditions~\eqref{eq:S1} and~\eqref{eq:S2} are equivalent if $N=1$ and $Q_1$ satisfies Slater's condition, i.e., there exists some $x_0\in\mathbb{R}^{n}$ such that $x_0^\top Q_1x_0>0$.
	\setcounter{thm}{\thelemma}\stepcounter{lemma}
\end{lem}

\begin{pf*}{Proof of Theorem~\ref{thm:observer-lmi}.}
	Let $M_{ij}$, $i,j\in\pazocal{N}$, be as defined in~\eqref{eq:lmi}. Suppose there exist symmetric matrices $P_i\in\mathbb{S}^{n_x}$, matrices $X_i\in\mathbb{R}^{n_x\times n_x}$, $Y_{ij}\in\mathbb{R}^{n_x\times n_y}$, $Z_{ij}\in\mathbb{R}^{n_\phi\times n_y}$, $i,j\in\pazocal{N}$, and scalars $\tau_{ij},\kappa_v,\kappa_w,\kappa_\psi\in\mathbb{R}_{>0}$, $i,j\in\pazocal{N}$, such that $M_{ij}\succ 0$ for all $i,j\in\pazocal{N}$. First, we note that it follows immediately from~\eqref{eq:lmi} that $P_i\succ I\succ 0$ for all $i\in\pazocal{N}$. Next, we show that the matrices $X_i$, $i\in\pazocal{N}$, are non-singular. To this end, let, for any $i\in\pazocal{N}$, $u\in\mathbb{R}^{n_x}$ be such that $X_i^\top u=0$, then, by~\eqref{eq:lmi}, $u^\top(\operatorname{He}(X_iE_j)-P_j)u=-u^\top P_ju{\geqslant} 0$ for all $j\in\pazocal{N}$. It follows, since $P(p)\succ 0$ for all $p\in\mathbb{P}$, that $u=0$ and, thus, the matrices $X_i$, $i\in\pazocal{N}$, are non-singular. Using the fact that, by Assumption~\ref{ass:parameter}, $\xi_i(p)\in\mathbb{R}_{\geqslant 0}$, $i\in\pazocal{N}$, and $\sum_{i\in\pazocal{N}}\xi_i(p)=1$ for all $p\in\pazocal{N}$, we have $P(p)\coloneqq\sum_{i\in\pazocal{N}}\xi_i(p)P_i\succ 0$, for all $p\in\mathbb{P}$.

We now show that a poly-quadratic ISS-Lyapunov function can be constructed.	Let $\Gamma\coloneqq\operatorname{diag}(\kappa_vI_{n_v},\kappa_wI_{n_w},\kappa_{\psi}I_{n_x})$. Using Assumption~\ref{ass:parameter}, we have, for all $i\in\pazocal{N}$ and $\hat{p}_+\in\mathbb{P}$,
	\begin{align}
		&0\prec\sum_{\mathclap{j\in\pazocal{N}}}\xi_j(\hat{p}_+)M_{ij} = \\
		&\left[\begin{smallarray}{@{}ccc;{2pt/2pt}c@{}}
			\operatorname{He}(X_iE(\hat{p}_+)) - P(\hat{p}_+) & \star & \star & \\
			(X_iA_i-Y_i(\hat{p}_+)C)^\top & P_i-I & \star & \,\star\\
			\vphantom{{X_i}_{\hat{\top}}}-(X_iG_i)^\top & \Lambda(\tau_i(\hat{p}_+)H-Z_i(\hat{p}_+)C) & 2\tau_i(\hat{p}_+)I & \\\hdashline[2pt/2pt]
			\vphantom{X_i^{\hat{\top}}}-(X_iF_i)^\top & 0 & 0 & \\
			(Y_i(\hat{p}_+)D)^\top & 0 & (\Lambda Z_i(\hat{p}_+)D)^\top & \,\Gamma\\
			X_i^\top & 0 & 0 &
		\end{smallarray}\right],\nonumber
	\end{align}
	where $\tau_i(\hat{p}_+) \coloneqq \sum_{{j\in\pazocal{N}}}\xi_j(\hat{p}_+)\tau_{ij}$, $Y_i(\hat{p}_+)\coloneqq \sum_{{j\in\pazocal{N}}}\xi_j(\hat{p}_+)Y_{ij}$ and $Z_i(\hat{p}_+)\coloneqq \sum_{{j\in\pazocal{N}}}\xi_j(\hat{p}_+)Z_{ij}$, $\hat{p}_+\in\mathbb{P}$ and $i\in\pazocal{N}$. It follows that $X_iE(p)+E^\top(p)X_i^\top - P(p) \succ 0$ for all $i\in\pazocal{N}$ and $p\in\mathbb{P}$. Therefore, $(P(p)-X_iE(p)) P^{-1}(p)(P(p)-X_iE(p))^\top\succcurlyeq 0$ for all $i\in\pazocal{N}$ and $p\in\mathbb{P}$ and, hence, for all $i\in\pazocal{N}$ and $p\in\mathbb{P}$, 
	\begin{equation}
		X_iE(p)P^{-1}(p)E^\top(p)X^\top_i\succcurlyeq \operatorname{He}(X_iE(p))-P(p)\succ 0.
		\label{eq:ineq}
	\end{equation} 	
	Using~\eqref{eq:ineq} followed by a congruence transformation with the non-singular matrix $T_i=\operatorname{diag}(X^{-\top}_i,I_{n_x},-I_{n_\phi},I_{n_x+n_v+n_w})$, $i\in\pazocal{N}$, $\hat{p}_+\in\mathbb{P}$, we have, for all $i\in\pazocal{N}$ and $\hat{p}_+\in\mathbb{P}$, that
	\begin{equation}
		\left[\begin{smallarray}{ccc;{2pt/2pt}c}
				E(\hat{p}_+)P^{-1}(\hat{p}_+)E^{\top}(\hat{p}_+) & \star & \star & \\
				\pazocal{A}^\top_i(\hat{p}_+) & P_i-I & \star & \,\star\\
				\vphantom{{X_i}_{\hat{\top}}}G_i^\top & -\Lambda\pazocal{H}_i(\hat{p}_+) & 2\tau_i(\hat{p}_+)I & \\\hdashline[2pt/2pt]
				\vphantom{F_i^{\hat{\top}}}-F_i^\top & 0 & 0 & \\
				(L_i(\hat{p}_+)D)^\top & 0 & -(\Lambda Z_i(\hat{p}_+)D)^\top & \,\Gamma\\
				I & 0 & 0 &
			\end{smallarray}\right] \succ 0,
		\label{eq:cond}
	\end{equation}
	where $\pazocal{A}_i(\hat{p}_+)\coloneqq A_i-L_i(\hat{p}_+)C$, $\pazocal{H}_i(\hat{p}_+)\coloneqq \tau_i(\hat{p}_+)H-Z_i(\hat{p}_+)C$, $L_i(\hat{p}_+)\coloneqq \sum_{j\in\pazocal{N}} \xi_j(\hat{p}_+)L_{ij}$ with $L_{ij} \coloneqq X_i^{-1}Y_{ij}$, $i,j\in\pazocal{N}$ and $\hat{p}_+\in\mathbb{P}$. Since~\eqref{eq:cond} is affine in all terms depending on $i$, we multiply by $\xi_i(\hat{p})$ and sum over all $i\in\pazocal{N}$ to find, for all $\hat{p}_+,\hat{p}\in\mathbb{P}$,
	\begin{equation}
		\left[\begin{smallarray}{ccc;{2pt/2pt}c}
				E(\hat{p}_+)P^{-1}(\hat{p}_+)E^\top(\hat{p}_+) & \star & \star & \\
				\pazocal{A}^\top(\hat{p}_+,\hat{p}) & P(\hat{p})-I & \star & \,\star\\
				\vphantom{{X_i}_{\hat{\top}}}G^\top(\hat{p}) & -\Lambda\pazocal{H}(\hat{p}_+,\hat{p}) & 2\tau(\hat{p}_+,\hat{p})I & \\\hdashline[2pt/2pt]
				\vphantom{F^{\hat{\top}}}-F^\top(\hat{p}) & 0 & 0 & \\
				(L(\hat{p}_+,\hat{p})D)^\top & 0 & -(\Lambda Z(\hat{p}_+,\hat{p})D)^\top & \,\Gamma\\
				I & 0 & 0 & 
			\end{smallarray}\right] \succ 0,
	\end{equation}
	with $\pazocal{A}(\hat{p}_+,\hat{p})\coloneqq A(\hat{p})-L(\hat{p}_+,\hat{p})C$, $\pazocal{H}(\hat{p}_+,\hat{p})\coloneqq \tau(\hat{p}_+,\hat{p})H-K(\hat{p}_+,\hat{p})C$, $\hat{p}_+,\hat{p}\in\mathbb{P}$, and $L$, $Z$ and $\tau$ as in~\eqref{eq:polytopic-L} and~\eqref{eq:K}-\eqref{eq:Z}. Using the Schur complement and the fact that $Z(\hat{p}_+,\hat{p})=\tau(\hat{p}_+,\hat{p})K(\hat{p}_+,\hat{p})$, we find, for all $\hat{p}_+,\hat{p}\in\mathbb{P}$,
	\begin{align}
		&\begin{bmatrix}\begin{smallmatrix}
			P(\hat{p})- I & \star & \star & \star & \star\\
			0 & 0 & \star & \star & \star\\
			0 & 0 & \kappa_vI & \star & \star\\
			0 & 0 & 0 & \kappa_wI & \star\\
			0 & 0 & 0 & 0 & \kappa_\psi I
		\end{smallmatrix}\end{bmatrix}-\label{eq:schur}\\
		&\quad \pazocal{B}^\top(\hat{p}_+,\hat{p}) E^{-\top}(\hat{p}_+)P(\hat{p}_+)E^{-1}(\hat{p}_+)\pazocal{B}(\hat{p}_+,\hat{p})-\nonumber\\
		&\quad\tau(\hat{p}_+,\hat{p})\begin{bmatrix}\begin{smallmatrix}
			0 & \star & \star & \star & \star\\
			\Lambda(H-K(\hat{p}_+,\hat{p})C) & -2 I & \star & \star & \star\\
			0 & 0 & 0 & \star & \star\\
			0 & (\Lambda K(\hat{p}_+,\hat{p})D)^\top & 0 & 0 & \star\\
			0 & 0 & 0 & 0 & 0
		\end{smallmatrix}\end{bmatrix}\succ 0,\nonumber
	\end{align}
	where, for $p_+,p\in\mathbb{P}$,
	\begin{equation}
		\pazocal{B}(p_+,p)\coloneqq \begin{bmatrix}\begin{smallmatrix}
			\pazocal{A}(\hat{p}_+,\hat{p}) & G(\hat{p}) & -F(\hat{p}) & L(\hat{p}_+,\hat{p})D & I
		\end{smallmatrix}\end{bmatrix}.
	\end{equation}
	Using pointwise application of Lemma~\ref{lem:S-proc},~\eqref{eq:schur} implies that, for all $x_+,x,e\in\mathbb{R}^{n_x}$, $v\in\mathbb{R}^{n_v}$, $w\in\mathbb{R}^{n_w}$, $u\in\mathbb{R}^{n_u}$ and $\hat{p}_+,\hat{p},p_+,p\in\mathbb{P}$,
	\begin{align}
		e_+^\top P(\hat{p}_+)e_+ - e^\top P(\hat{p})e &\leqslant - \|e\|^2 + \label{eq:interm1}\\
		&\quad\kappa_v\|v\|^2+\kappa_w\|w\|^2+\kappa_{\psi}\|\psi\|^2,\nonumber
	\end{align}
	for all $\tilde{\phi}$ subject to $\|\tilde{\phi}\|^2-\tilde{\phi}^\top \Lambda ((H-K(\hat{p}_+,\hat{p})C)e+K(\hat{p}_+,\hat{p})Dw)\leqslant 0$, which holds for any $\tilde{\phi}$ satisfying~\eqref{eq:prop-slope-restr} with $\phi(y)-\phi(x)\leftarrow\tilde{\phi}$, $x\leftarrow Hx$ and $y\leftarrow(H-K(\hat{p}_+,\hat{p})C)e+K(\hat{p}_+,\hat{p})Dw+Hx$, and $e_+$ subject to $E(\hat{p}_+)e_+ = \pazocal{B}(\hat{p}_+,\hat{p})q$ with 
	\begin{equation}
		q\coloneqq (e,\tilde{\phi},v,w,\psi),
		\label{eq:q}
	\end{equation}
	$\tilde{\phi}=\tilde{\phi}(\hat{p}_+,\hat{p},e,x,w)$ and $\psi=\psi(\hat{p}_+,\hat{p},p_+,p,x_+,x,u,v)$. In other words, $V(e,p) = e^\top P(p)e$ satisfies~\eqref{eq:descent-cond}. $V$ also satisfies~\eqref{eq:lyap-bounds}, since $P(\hat{p})\succ I$, and, from the bottom-right of~\eqref{eq:schur}, $\kappa_\psi E^\top(\hat{p}_+)E(\hat{p}_+) \succ P(\hat{p}_+)$ for all $\hat{p}_+\in\mathbb{P}$, which implies $I\prec P(p)\prec \kappa_\psi \underline{\sigma}^2I$, with $\underline{\sigma}$ as in Theorem~\ref{thm:observer-lmi}, for all $p\in\mathbb{P}$. Thus, $\kappa_\psi > \underline{\sigma}^{-2}$ and $V$ is a poly-quadratic ISS-Lyapunov function and the system is (absolutely) poly-quadratically ISS.
	
	Finally, we show ISS and derive the ISS-gains. Denote $V_k=V(e_k,\hat{p}_k)$ and apply~\eqref{eq:lyap-bounds} and~\eqref{eq:descent-cond} to find $V_{k+1} \leqslant (1-\frac{1}{\kappa_\psi\underline{\sigma}^2})V_k + \kappa_v\|v_k\|^2 + \kappa_w\|w_k\|^2 + \kappa_\psi\|\psi_k\|^2$. Applying repetitively yields	$V_{k} \leqslant \rho^{k}V_0 + \kappa_\psi\underline{\sigma} ^2(\kappa_v\|v\|^2_{\infty,k} + \kappa_w\|w\|^2_{\infty,k} + \kappa_{\psi}\|\psi\|^2_{\infty,k})$, with $\rho\coloneqq(1-(\kappa_{\psi}\underline{\sigma}^2)^{-1})\in(0,1)$ since $1 < \kappa_{\psi}\underline{\sigma}^2$. Using~\eqref{eq:lyap-bounds} and $\sqrt{a+b}\leqslant \sqrt{a}+\sqrt{b}$ for any $a,b\in\mathbb{R}_{\geqslant 0}$, we have ISS with respect to $v$, $w$ and $\psi$ with $\beta\in\pazocal{KL}$ and $\gamma_v,\gamma_w,\gamma_\psi\in\pazocal{K}$ as in Theorem~\ref{thm:observer-lmi}.
\end{pf*}

\begin{pf*}{Proof of Theorem~\ref{thm:suffnec}.}
	We prove \underline{sufficiency} following similar steps as in the proof of Theorem~\ref{thm:observer-lmi}. Let $N_{ij}$, $i,j\in\pazocal{N}$, be as defined in~\eqref{eq:nec-cond}, and suppose there exist matrices $P_i\in\mathbb{S}^{n_x}$, $X_{ij}\in\mathbb{R}^{n_x\times n_x}$, $Y_{ij}\in\mathbb{R}^{n_x\times n_y}$, $Z_{ij}\in\mathbb{R}^{n_\phi\times n_y}$ and scalars $\tau_{ij},\kappa_v,\kappa_w,\kappa_{\psi}\in\mathbb{R}_{>0}$, $i,j\in\pazocal{N}$, such that $N_{ij}\succ 0$ for all $i,j\in\pazocal{N}$. First, it follows immediately from~\eqref{eq:nec-cond} that $P_i\succ I\succ 0$, for all $i\in\pazocal{N}$, and, hence, $P(p)\succ 0$ for all $p\in\mathbb{P}$. Following the same reasoning as in the proof of Theorem~\ref{thm:observer-lmi}, we can show that the matrices $X_{ij}$, $i,j\in\pazocal{N}$, are non-singular and, thus, $X_i(\hat{p}_+)$, $i\in\pazocal{N}$, as in~\eqref{eq:XY} are also non-singular for all $\hat{p}_+\in\mathbb{P}$. 
	
	Next, we show that a poly-quadratic ISS-Lyapunov function can be constructed. Let $\Gamma\coloneqq \operatorname{diag}(\kappa_vI_{n_v},\kappa_wI_{n_w},\kappa_{\psi}I_{n_x})$. Using Assumption~\ref{ass:parameter}, $\xi_i(p)\in\mathbb{R}_{\geqslant 0}$, $i\in\pazocal{N}$, and $\sum_{i\in\pazocal{N}}\xi_i(p)=1$ for all $p\in\mathbb{P}$, we have, for all $i\in\pazocal{N}$ and $\hat{p}_+\in\mathbb{P}$,
	\begin{align}
		&0\prec\sum_{j\in\pazocal{N}}\xi_j(\hat{p}_+)N_{ij} = \\
		&\left[\begin{smallarray}{ccc;{2pt/2pt}c}
			\operatorname{He}(X_i(\hat{p}_+)\bar{E})-P(\hat{p}_+) & \star & \star & \\
			(X_i(\hat{p}_+)A_i-Y_i(\hat{p}_+)C)^\top & P_i-I & \star & \,\star\\
			\vphantom{{X_i}_{\hat{\top}}}-(X_i(\hat{p}_+)G_i)^\top & \Lambda(\tau_i(\hat{p}_+)H-Z_i(\hat{p}_+)C) & 2\tau_i(\hat{p}_+)I & \\\hdashline[2pt/2pt]
			\vphantom{X_i^{\hat{\top}}}-(X_i(\hat{p}_+)F_i)^\top & 0 & 0 & \\
			(Y_i(\hat{p}_+)D)^\top & 0 & (\Lambda Z_i(\hat{p}_+)D)^\top & \,\Gamma\\
			X_i^\top(\hat{p}_+) & 0 & 0 & 
		\end{smallarray}\right],\nonumber
	\end{align}
	where $Z_i(\hat{p}_+)\coloneqq\sum_{j\in\pazocal{N}}\xi_j(\hat{p}_+)Z_{ij}$, $\hat{p}_+\in\mathbb{P}$, $i\in\pazocal{N}$, and $X_i(\hat{p}_+)$ and $Y_i(\hat{p}_+)$ are as defined in~\eqref{eq:XY}. It follows that $P(p)\coloneqq \sum_{i\in\pazocal{N}}\xi_i(p)P_i\succ I \succ 0$ and $(P(p)-X_i\bar{E})P^{-1}(p)(P(p)-X_i\bar{E})^\top\succcurlyeq 0$ for all $i\in\pazocal{N}$ and $p\in\mathbb{P}$. Thus, it holds, for all $p\in\mathbb{P}$ and $i\in\pazocal{N}$, that
	\begin{equation}
		X_i(p)\bar{E}P^{-1}(p)\bar{E}^\top X_i^\top(p) \succcurlyeq \operatorname{He}(X_i(p)\bar{E})-P(p)\succ 0.
		\label{eq:ineq2}
	\end{equation}
	Using~\eqref{eq:ineq2} and congruence transformations with the non-singular matrices $T_i(\hat{p}_+)=\operatorname{diag}(X_i^{-\top}(\hat{p}_+),I_{n_x},-I_{n_\phi},I_{n_x+n_v+n_w})$, $i\in\pazocal{N}$, $\hat{p}_+\in\mathbb{P}$, we obtain, for all $\hat{p}_+\in\mathbb{P}$ and $i\in\pazocal{N}$,
	\begin{equation}
		\left[\begin{smallarray}{ccc;{2pt/2pt}c}
			\bar{E}P^{-1}(\hat{p}_+)\bar{E}^\top & \star & \star & \\
			\pazocal{A}_i^\top(\hat{p}_+) & P_i-I & \star & \,\star\\
			\vphantom{{X_i}_{\hat{\top}}}G_i^\top & -\Lambda\pazocal{H}_i(\hat{p}_+) & 2\tau_i(\hat{p}_+)I & \\\hdashline[2pt/2pt]
			\vphantom{F_i^{\hat{\top}}}-F_i^\top & 0 & 0 & \\
			(L_i(\hat{p}_+)D)^\top & 0 & -(\Lambda Z_i(\hat{p}_+)D)^\top & \,\Gamma\\
			I & 0 & 0 & 
		\end{smallarray}\right]\succ 0.
		\label{eq:affine}
	\end{equation}
	where $\pazocal{A}_i(\hat{p}_+)\coloneqq A_i-L_i(\hat{p}_+)C$, $\pazocal{H}_i(\hat{p}_+)\coloneqq \tau_i(\hat{p}_+)H-Z_i(\hat{p}_+)C$, $L_i(\hat{p}_+) \coloneqq X_i^{-1}(\hat{p}_+)Y_i(\hat{p}_+)$, $\hat{p}_+\in\mathbb{P}$, $i\in\pazocal{N}$ and $\hat{p}_+\in\mathbb{P}$. Since~\eqref{eq:affine} is affine in all terms depending on $i$, multiplying by $\xi_i(\hat{p})$ and summing over all $i\in\pazocal{N}$ yields, for all $\hat{p}_+,\hat{p}\in\mathbb{P}$,
	\begin{equation}
		\left[\begin{smallarray}{ccc;{2pt/2pt}c}
			\bar{E}P^{-1}(\hat{p}_+)\bar{E}^\top & \star & \star & \\
			\pazocal{A}^\top(\hat{p}_+,\hat{p}) & P(\hat{p})-I & \star & \,\star\\
			\vphantom{{X_i}_{\hat{\top}}}G^\top(\hat{p}) & -\Lambda\pazocal{H}(\hat{p}_+,\hat{p}) & 2\tau(\hat{p}_+,\hat{p})I & \\\hdashline[2pt/2pt]
			\vphantom{F^{\hat{\top}}}-F^\top(\hat{p}) & 0 & 0 & \\
			(L(\hat{p}_+,\hat{p})D)^\top & 0 & -(\Lambda Z(\hat{p}_+,\hat{p})D)^\top & \,\Gamma\\
			I & 0 & 0 & 
		\end{smallarray}\right]\succ 0,
	\end{equation}
	with $\pazocal{A}(\hat{p}_+,\hat{p})\coloneqq A(\hat{p})-L(\hat{p}_+,\hat{p})C$, $\pazocal{H}(\hat{p}_+,\hat{p})\coloneqq \tau(\hat{p}_+,\hat{p})H-Z(\hat{p}_+,\hat{p})C$, $\hat{p}_+,\hat{p}\in\mathbb{P}$, and $L$, $Z$ and $\tau$ as in~\eqref{eq:L}-\eqref{eq:XY} and~\eqref{eq:K}-\eqref{eq:Z}. Applying the Schur complement and the fact that $Z(\hat{p}_+,\hat{p})=\tau(\hat{p}_+,\hat{p})K(\hat{p}_+,\hat{p})$, we find, for all $\hat{p}_+,\hat{p}\in\mathbb{P}$,
	\begin{align}
			&\begin{bmatrix}\begin{smallmatrix}
				P(\hat{p})-I & \star & \star & \star & \star\\
				0 & 0 & \star & \star & \star\\
				0 & 0 & \kappa_vI & \star & \star\\
				0 & 0 & 0 & \kappa_wI & \star\\
				0 & 0 & 0 & 0 & \kappa_{\psi}I
			\end{smallmatrix}\end{bmatrix} - 
		\label{eq:pre-S-proc}\\
		&\quad\pazocal{B}^\top(\hat{p}_+,\hat{p})\bar{E}^{-\top}P(\hat{p}_+)\bar{E}^{-1}\pazocal{B}(\hat{p}_+,\hat{p})-\nonumber\\
			&\quad \tau(\hat{p}_+,\hat{p})\begin{bmatrix}\begin{smallmatrix}
				0 & \star & \star & \star & \star\\
				\Lambda (H-K(\hat{p}_+,\hat{p})C) & -2I & \star & \star & \star\\
				0 & 0 & 0 & \star & \star\\
				0 & (\Lambda K(\hat{p}_+,\hat{p})D)^\top & 0 & 0 & \star\\
				0 & 0 & 0 & 0 & 0
			\end{smallmatrix}\end{bmatrix} \succ 0,\nonumber
	\end{align}
	where, for $p_+,p\in\mathbb{P}$,
	\begin{equation}
		\pazocal{B}(\hat{p}_+,\hat{p})\coloneqq\begin{bmatrix}\begin{smallmatrix}
			\pazocal{A}(\hat{p}_+,\hat{p}) & G(\hat{p}) & -F(\hat{p}) & L(\hat{p}_+,\hat{p})D & I
		\end{smallmatrix}\end{bmatrix}.
		\label{eq:Bcal}
	\end{equation}
	Using pointwise application of Lemma~\ref{lem:S-proc}, it follows that, for all $x_+,x,e\in\mathbb{R}^{n_x}$, $v\in\mathbb{R}^{n_v}$, $w\in\mathbb{R}^{n_w}$, $u\in\mathbb{R}^{n_u}$ and $\hat{p}_+,\hat{p},p_+,p\in\mathbb{P}$,~\eqref{eq:interm1} holds for all $\tilde{\phi}$ subject to $\|\tilde{\phi}\|^2-\tilde{\phi}^\top\Lambda((H-K(\hat{p}_+,\hat{p})C)e+K(\hat{p}_+,\hat{p})Dw)\leqslant 0$, which holds for any $\tilde{\phi}$ satisfying~\eqref{eq:prop-slope-restr} with $\phi(y)-\phi(x)\leftarrow\tilde{\phi}$, $x\leftarrow Hx$ and $y\leftarrow(H-K(\hat{p}_+,\hat{p})C)e+K(\hat{p}_+,\hat{p})Dw+Hx$, and $e_+$ subject to $\bar{E}e_+=\pazocal{B}(\hat{p}_+,\hat{p})q$ with $q$ as in~\eqref{eq:q}, $\tilde{\phi}=\tilde{\phi}(\hat{p}_+,\hat{p},e,x,w)$ and $\psi=\psi(\hat{p}_+,\hat{p},p_+,p,x_+,x,u,v)$. From the bottom-right of~\eqref{eq:pre-S-proc}, we have that $P(p)\prec \kappa_\psi \underline{\sigma}^2I$ with $\underline{\sigma}$ as in Theorem~\ref{thm:suffnec}, for all $p\in\mathbb{P}$. Thus, $V$ is a poly-quadratic ISS-Lyapunov function (as shown in more detail in the proof of Theorem~\ref{thm:observer-lmi}) and the error system is (absolutely) poly-quadratically ISS. The precise ISS-gains can also be derived as in Theorem~\ref{thm:observer-lmi}.
	
	To show \underline{necessity}, suppose there exist matrix-valued functions $L\colon\mathbb{P}\times\mathbb{P}\rightarrow\mathbb{R}^{n_x\times n_y}$ and $K\colon\mathbb{P}\times\mathbb{P}\rightarrow\mathbb{R}^{n_\phi\times n_y}$ for which the observer~\eqref{eq:observer} renders the error system~\eqref{eq:error-system} absolutely poly-quadratically ISS with respect to $v$, $w$ and $\psi$. Then, by Definition~\ref{def:poly-ISS}, there exist symmetric matrices $\bar{P}_i$, $i\in\pazocal{N}$, such that $\bar{V}(e,\hat{p})=e^\top\bar{P}(\hat{p})e$ with $\bar{P}(p)\coloneqq \sum_{i\in\pazocal{N}}\xi_i(p)\bar{P}_i$, $p\in\mathbb{P}$, satisfies $0\prec I\preccurlyeq \bar{P}(\hat{p})\preccurlyeq \bar{a}I$, for all $\hat{p}\in\mathbb{P}$ and for some $\bar{a}\in\mathbb{R}_{>0}$. Since $\{\bm{e}_i\}_{i\in\pazocal{N}}\subset\xi(\mathbb{P})$, it follows that $0\prec I\preccurlyeq \bar{P}_i\preccurlyeq\bar{a}I$ for all $i\in\pazocal{N}$. Then, from~\eqref{eq:descent-cond}, we have, for all $\hat{p}_+,\hat{p}\in\mathbb{P}$, 
	\begin{align}
		&q^\top\pazocal{B}^\top(\hat{p}_+,\hat{p})\bar{E}^{-\top}\bar{P}(\hat{p}_+)\bar{E}^{-1}\pazocal{B}(\hat{p}_+,\hat{p})q - e^\top\bar{P}(\hat{p})e\leqslant\nonumber \\
		&\quad-\|e\|^2+\bar{\kappa}_v\|v\|^2+\bar{\kappa}_w\|w\|^2+\bar{\kappa}_{\psi}\|\psi\|^2,
		\label{eq:star}
	\end{align}
	with $\pazocal{B}$ as defined in~\eqref{eq:Bcal}, for some $\bar{\kappa}_v,\bar{\kappa}_w,\bar{\kappa}_{\psi}\in\mathbb{R}_{>0}$ and for all $q\in\mathbb{R}^{n_q}$, with $n_q=2n_x+n_u+n_v+n_w$, subject to
	\begin{equation}
		q^\top \underbrace{\begin{bmatrix}\begin{smallmatrix}
			0 & \star & \star & \star & \star\\
			\Lambda(H-K(\hat{p}_+,\hat{p})C) & -2I & \star & \star & \star\\
			0 & 0 & 0 & \star & \star\\
			0 & (\Lambda K(\hat{p}_+,\hat{p})D)^\top & 0 & 0 & \star\\
			0 & 0 & 0 & 0 & 0
		\end{smallmatrix}\end{bmatrix}}_{\eqqcolon \pazocal{W}(\hat{p}_+,\hat{p})}q\geqslant 0.
		\label{eq:S-proc}
	\end{equation}
	Note that~\eqref{eq:S-proc} is Assumption~\ref{ass:nonlinearity-condition} with $\phi(y)-\phi(x)\leftarrow\tilde{\phi}$, $x\leftarrow Hx$ and $y\leftarrow(H-K(\hat{p}_+,\hat{p})C)e+K(\hat{p}_+,\hat{p})Dw+Hx$ expressed in terms of $q$. We proceed to show necessity by taking the following steps: \begin{enumerate}[label=\alph*)]
		\item We construct $P\colon\mathbb{P}\rightarrow\mathbb{S}^{n_x}_{\succ 0}$ and $\kappa_v,\kappa_w,\kappa_\psi\in\mathbb{R}_{>0}$ such that, for all $\hat{p}_+,\hat{p}\in\mathbb{P}$, $q^\top\pazocal{V}(\hat{p}_+,\hat{p})q>0$ for all $q\in\mathbb{R}^{n_q}\setminus\{0\}$ when $q^\top\pazocal{W}(\hat{p}_+,\hat{p})q\geqslant 0$, where
		\begin{align}
			&\pazocal{V}(\hat{p}_+,\hat{p})\coloneqq \operatorname{diag}(P(\hat{p})-I,0,\kappa_vI,\kappa_wI,\kappa_{\psi}I)-\nonumber\\
			&\qquad\pazocal{B}^\top(\hat{p}_+,\hat{p})\bar{E}^{-\top}P(\hat{p}_+)\bar{E}^{-1}\pazocal{B}(\hat{p}_+,\hat{p}),
			\label{eq:Vcal}
	\end{align}
	for $\hat{p}_+,\hat{p}\in\mathbb{P}$. We achieve this using the following lemma, for which a proof is provided later in the Appendix.
	\setcounter{thm}{\thelemma}\stepcounter{lemma}
	\begin{lem}\label{lem:shiftVcal}
		Consider $\bar{P}\colon\mathbb{P}\rightarrow\mathbb{S}^{n_x}_{\succ 0}$ and $\bar{\kappa}_v,\bar{\kappa}_w,\bar{\kappa}_{\psi}\in\mathbb{R}_{>0}$ such that, for all $\hat{p}_+,\hat{p}\in\mathbb{P}$,~\eqref{eq:star} holds for all $q\in\mathbb{R}^{n_q}$ subject to~\eqref{eq:S-proc}. For any $\epsilon\in(0,1)$, we have, for all $\hat{p}_+,\hat{p}\in\mathbb{P}$,
		\begin{equation}
			q^\top\bar{\pazocal{V}}(\hat{p}_+,\hat{p})q>0,
			\label{eq:shifted}
		\end{equation}
		with $\bar{\pazocal{V}}(\hat{p}_+,\hat{p})\coloneqq \operatorname{diag}(\bar{P}(\hat{p})-(1-\epsilon)I,0,(\bar{\kappa}_v+\epsilon)I,(\bar{\kappa}_w+\epsilon)I,(\bar{\kappa}_{\psi}+\epsilon)I)-\pazocal{B}^\top(\hat{p}_+,\hat{p})\bar{E}^{-\top}\bar{P}(\hat{p}_+)\bar{E}^{-1}\pazocal{B}(\hat{p}_+,\hat{p})$, for all $q\in\mathbb{R}^{n_q}\setminus\{0\}$ subject to~\eqref{eq:S-proc}.
	\end{lem}
		\item Next, we use the following $\pazocal{S}$-procedure-like result, for which a proof is provided later in the Appendix, to obtain a single matrix inequality that is equivalent to $q^\top \pazocal{V}(\hat{p}_+,\hat{p})q>0$ for all $q\in\mathbb{R}^{n_q}\setminus\{0\}$ when $q^\top\pazocal{W}(\hat{p}_+,\hat{p})q\geqslant 0$ with $\pazocal{V}$ and $\pazocal{W}$ as in~\eqref{eq:Vcal} and~\eqref{eq:S-proc}, respectively.
		\setcounter{thm}{\thelemma}\stepcounter{lemma}
		\begin{lem}\label{lem:tau}
			Consider $\pazocal{V}\colon\mathbb{P}\times\mathbb{P}\rightarrow\mathbb{S}^{n_q}$ and $\pazocal{W}\colon\mathbb{P}\times\mathbb{P}\rightarrow\mathbb{S}^{n_q}$ as in~\eqref{eq:Vcal} and~\eqref{eq:S-proc}, respectively. There exists a function $\tau\colon\mathbb{P}\times\mathbb{P}\rightarrow\mathbb{R}_{>0}$ such that, for all $\hat{p}_+,\hat{p}\in\mathbb{P}$,
			\begin{equation}
				\pazocal{V}(\hat{p}_+,\hat{p})-\tau(\hat{p}_+,\hat{p})\pazocal{W}(\hat{p}_+,\hat{p})\succ 0,
			\label{eq:VtauW}
			\end{equation}
			if and only if, for all $\hat{p}_+,\hat{p}\in\mathbb{P}$, it holds that
			\begin{align}
				&q^\top\pazocal{V}(\hat{p}_+,\hat{p})q>0,\text{ for all }q\in\mathbb{R}^{n_q}\setminus\{0\}\nonumber\\
				&\qquad\text{such that }q^\top\pazocal{W}(\hat{p}_+,\hat{p})q\geqslant 0.
			\end{align}
		\end{lem}
		At first glance, Lemma~\ref{lem:tau} might appear to simply be a pointwise application of Lemma~\ref{lem:S-proc}. However, although necessity follows immediately from Lemma~\ref{lem:S-proc}, sufficiency is not immediately clear since $\pazocal{W}(\hat{p}_+,\hat{p})$ does not necessarily satisfy Slater's condition for all $\hat{p}_+,\hat{p}\in\mathbb{P}$. In fact, we distinguish, for each $\hat{p}_+,\hat{p}\in\mathbb{P}$, two cases: \begin{enumerate*}[label=(\roman*)]
			\item $(H-K(\hat{p}_+,\hat{p})C)\neq 0$ or $K(\hat{p}_+,\hat{p})D\neq 0$, or
			\item $(H-K(\hat{p}_+,\hat{p})C)=0$ and $K(\hat{p}_+,\hat{p})D=0$.		
		\end{enumerate*}
		In the proof of Lemma~\ref{lem:tau}, we apply pointwise for each $\hat{p}_+,\hat{p}\in\mathbb{P}$ either \begin{enumerate*}[label=(\roman*)]
			\item Lemma~\ref{lem:S-proc} since Slater's condition holds in this case, or
			\item Lemma~\ref{lem:finsler},
		\end{enumerate*}
		to construct $\tau\colon\mathbb{P}\times\mathbb{P}\rightarrow\mathbb{R}_{>0}$ satisfying~\eqref{eq:VtauW}. 
		\item Finally, we exploit the fact that $\{\bm{e}_i\}_{i\in\pazocal{N}}\subset\xi(\mathbb{P})$ to evaluate the matrix inequality obtained in the previous step, which depends on $\hat{p}_+,\hat{p}\in\mathbb{P}$, to obtain a finite set of conditions corresponding to each of the vertices. We then apply Lemma~\ref{lem:proj} to arrive at the conditions in Theorem~\ref{thm:suffnec} and, hence, complete the necessity part of the proof. 
	\end{enumerate}
	Let us now proceed along the steps outlined above. As for step a), we directly obtain~\eqref{eq:Vcal} using Lemma~\ref{lem:shiftVcal} by taking $P(p)\coloneqq\bar{P}(p)/(1-\epsilon)$, $p\in\mathbb{P}$, $\kappa_v\coloneqq (\bar{\kappa}_v+\epsilon)/(1-\epsilon)$, $\kappa_w\coloneqq (\bar{\kappa}_w+\epsilon)/(1-\epsilon)$ and $\kappa_{\psi}\coloneqq (\bar{\kappa}_{\psi}+\epsilon)/(1-\epsilon)$. Next, we can directly apply Lemma~\ref{lem:tau} to conclude that there exists a function $\tau\colon\mathbb{P}\times\mathbb{P}\rightarrow\mathbb{R}_{>0}$ such that~\eqref{eq:VtauW} holds for all $\hat{p}_+,\hat{p}\in\mathbb{P}$. 
	
	We now proceed with step c) from our earlier outline. Since $\{\bm{e}_i\}_{i\in\pazocal{N}}\subset\xi(\mathbb{P})$, there exist $\nu_i\in\mathbb{P}$ for which $\bm{e}_i=\xi(\nu_i)$, $i\in\pazocal{N}$. Substitution in~\eqref{eq:VtauW}, yields, for all $i,j\in\pazocal{N}$,
	\begin{align}
		&0\prec\underbrace{\begin{bmatrix}\begin{smallmatrix}
			P_i-I & \star & \star & \star & \star\\
			-\Lambda(\tau_{ij}H-Z_{ij}C) & 2\tau_{ij}  & \star & \star & \star\\
			0 & 0 & \kappa_vI & \star & \star\\
			0 & -(\Lambda Z_{ij}D)^\top & 0 & \kappa_wI & \star \\
			0 & 0 & 0 & 0 & \kappa_\psi I 
		\end{smallmatrix}\end{bmatrix}}_{\eqqcolon \pazocal{Q}_{ij}} -\label{eq:temp}\\
		&\quad\begin{bmatrix}\begin{smallmatrix}
			(A_i-L_{ij}C)^\top\\
			G_i^\top\\
			-F_i^\top\\
			(L_{ij}D)^\top\\
			I
		\end{smallmatrix}\end{bmatrix}\bar{E}^{-\top}P_j\bar{E}^{-1}\underbrace{\begin{bmatrix}\begin{smallmatrix}
			(A_i-L_{ij}C) & G_i & -F_i & L_{ij}D & I
		\end{smallmatrix}\end{bmatrix}}_{\eqqcolon\pazocal{B}_{ij}},\nonumber		
	\end{align}
	where $\tau_{ij} = \tau(\nu_j,\nu_i)$, $Z_{ij} = \tau_{ij}K(\nu_j,\nu_i)$, $P_i=P(\nu_i)$ and $L_{ij}=L(\nu_{j},\nu_i)$, $i,j\in\pazocal{N}$. We express~\eqref{eq:temp} equivalently as, for all $i,j\in\pazocal{N}$,
	\begin{equation}
		\begin{bmatrix}\begin{smallmatrix}
			-\pazocal{B}_{ij}^\top \bar{E}^{-\top} & I
		\end{smallmatrix}\end{bmatrix}\begin{bmatrix}\begin{smallmatrix}
			-P_j & 0\\
			0 & \pazocal{Q}_{ij}
		\end{smallmatrix}\end{bmatrix}\begin{bmatrix}\begin{smallmatrix}
			-\bar{E}^{-1}\pazocal{B}_{ij}\\ I
		\end{smallmatrix}\end{bmatrix} \succ 0.
		\label{eq:proj1}
	\end{equation}
	It follows that $Q_{ij}\succ \pazocal{B}_{ij}^\top\bar{E}^{-\top}P_j\bar{E}_j^{-1}\pazocal{B}_{ij}\succcurlyeq 0$ and, thus, it holds that, for all $i,j\in\pazocal{N}$,
	\begin{equation}
		\begin{bmatrix}\begin{smallmatrix}
			0 & I
		\end{smallmatrix}\end{bmatrix}\begin{bmatrix}\begin{smallmatrix}
			-P_j & 0\\
			0 & \pazocal{Q}_{ij}
		\end{smallmatrix}\end{bmatrix}\begin{bmatrix}\begin{smallmatrix}
			0\\
			I
		\end{smallmatrix}\end{bmatrix} \succ 0.
		\label{eq:proj2}
	\end{equation}
	Applying Lemma~\ref{lem:proj} to~\eqref{eq:proj1} and~\eqref{eq:proj2}, there exist matrices $X_{ij}\in\mathbb{R}^{}$, $i,j\in\pazocal{N}$, such that, for all $i,j\in\pazocal{N}$, we have
	\begin{equation}
		\begin{bmatrix}\begin{smallmatrix}
			-P_j & 0\\
			0 & \pazocal{Q}_{ij}
		\end{smallmatrix}\end{bmatrix} + \begin{bmatrix}\begin{smallmatrix}
			I\\
			0
		\end{smallmatrix}\end{bmatrix}X_{ij}\begin{bmatrix}\begin{smallmatrix}
			\bar{E} & \pazocal{B}_{ij}
		\end{smallmatrix}\end{bmatrix} + \begin{bmatrix}\begin{smallmatrix}
			\bar{E}^\top\\
			\pazocal{B}_{ij}^\top
		\end{smallmatrix}\end{bmatrix}X_{ij}^\top\begin{bmatrix}\begin{smallmatrix}
			I & 0
		\end{smallmatrix}\end{bmatrix}\succ 0.
	\end{equation}
	Performing a congruence transformation with $\operatorname{diag}(I_{2n_x},-I_{n_\phi},I_{n_v+n_w+n_x})$, we have, for all $i,j\in\pazocal{N}$,
	\begin{equation}
		\begin{bmatrix}\begin{smallmatrix}
			\operatorname{He}(X_{ij}\bar{E})-P_j & \star & \star & \star & \star\\
			(X_{ij}A_i-Y_{ij}C)^\top & P_i-I & \star & \star & \star & \star\\
			-(X_{ij}G_i)^\top & \Lambda(\tau_{ij}H - Z_{ij}C) & 2\tau_{ij}I & \star & \star & \star\\
			-(X_{ij}F_i)^\top & 0 & 0 & \kappa_vI & \star & \star\\
			(Y_{ij}D)^\top & 0 & (\Lambda Z_{ij}D)^{\top} & 0 & \kappa_wI & \star\\
			X_{ij}^\top & 0 & 0 & 0 & 0 & \kappa_{\psi}I
		\end{smallmatrix}\end{bmatrix}\succ 0,
		\label{eq:postcongr}
	\end{equation}
	where $Y_{ij}=X_{ij}L_{ij}$, $i,j\in\pazocal{N}$. Note that~\eqref{eq:postcongr} implies that $\tau_{ij}>0$ for all $i,j\in\pazocal{N}$. Thus, we have completed the proof of necessity by construction of the matrices $P_i\in\mathbb{S}^{n_x}$, $X_{ij}\in\mathbb{R}^{n_x\times n_x}$, $Y_{ij}\in\mathbb{R}^{n_x\times n_y}$, $Z_{ij}\in\mathbb{R}^{n_\phi\times n_y}$ and scalars $\tau_{ij},\kappa_v,\kappa_w,\kappa_{\psi}\in\mathbb{R}_{>0}$, $i,j\in\pazocal{N}$ that satisfy~\eqref{eq:nec-cond} for all $i,j\in\pazocal{N}$.
\end{pf*}

\begin{pf*}{Proof of Lemma~\ref{lem:shiftVcal}.}
	Let $\bar{P}\colon\mathbb{P}\rightarrow\mathbb{S}^{n_x}_{\succ 0}$ and $\bar{\kappa}_v,\bar{\kappa}_w,\bar{\kappa}_\psi\in\mathbb{R}_{>0}$ be such that, for all $\hat{p}_+,\hat{p}\in\mathbb{P}$,~\eqref{eq:star} holds for all $q\in\mathbb{R}^{n_q}$ subject to~\eqref{eq:S-proc}. To show that, for any $\epsilon\in(0,1)$,~\eqref{eq:shifted} holds for all $q\in\mathbb{R}^{n_q}\setminus\{0\}$ subject to~\eqref{eq:S-proc}, suppose, to the contrary, that there exist some $\bar{q}\in\mathbb{R}^{n_q}\setminus\{0\}$ and some $\hat{p}_+,\hat{p}\in\mathbb{P}$ such that $\bar{q}^\top\pazocal{W}(\hat{p}_+,\hat{p})\bar{q}\geqslant 0$ and $\bar{q}^\top\bar{\pazocal{V}}(\hat{p}_+,\hat{p})\bar{q}\leqslant 0$. Since we know already, using~\eqref{eq:star} and $\epsilon>0$, that, for all $\hat{p}_+,\hat{p}\in\mathbb{P}$, $\bar{q}^\top\bar{\pazocal{V}}(\hat{p}_+,\hat{p})\bar{q}\geqslant 0$ for all $\bar{q}\in\mathbb{R}^{n_q}$ subject to $\bar{q}^\top\pazocal{W}(\hat{p}_+,\hat{p})\bar{q}\geqslant 0$, it must hold that $\bar{q}^\top\bar{\pazocal{V}}(\hat{p}_+,\hat{p})\bar{q}=0$ for some $\hat{p}_+,\hat{p}\in\mathbb{P}$. By inspection of~\eqref{eq:shifted} using~\eqref{eq:star} (particularly, how we added $\epsilon$ to specific diagonal blocks), we see that any $\bar{q}\neq 0$ for which $\bar{q}^\top\pazocal{W}(\hat{p}_+,\hat{p})\bar{q}\geqslant 0$ and $\bar{q}^\top\bar{\pazocal{V}}(\hat{p}_+,\hat{p})\bar{q}=0$, for some $\hat{p}_+,\hat{p}\in\mathbb{P}$, must satisfy
	\begin{equation}
		\bar{q}\in\operatorname{Im}\underbrace{\begin{bmatrix}\begin{smallmatrix}
			0\\
			I_{n_{\phi}}\\
			0\\
			0\\
			0
		\end{smallmatrix}\end{bmatrix}}_{\eqqcolon \pazocal{M}},
		\label{eq:M}
	\end{equation}
	and, hence, we can express $\bar{q}$ as $\bar{q}=\pazocal{M}\bar{v}$ for some $\bar{v}\in\mathbb{R}^{n_\phi}\setminus\{0\}$. Substitution in $\bar{q}^\top\pazocal{W}(\hat{p}_+,\hat{p})\bar{q}\geqslant 0$ shows that $0\leqslant \bar{v}^\top\pazocal{M}^\top\pazocal{W}(\hat{p}_+,\hat{p})\pazocal{M}\bar{v}=-2\|\bar{v}\|^2$, which implies that $\bar{v}=0$ and, hence, $\bar{q}=0$. This contradicts our initial statement in which $\bar{q}\neq 0$, which completes our proof that~\eqref{eq:shifted} holds for all $q\in\mathbb{R}^{n_q}\setminus\{0\}$ subject to~\eqref{eq:S-proc}.	
\end{pf*}

\begin{pf*}{Proof of Lemma~\ref{lem:tau}.}
	Since necessity readily follows from pointwise application of Lemma~\ref{lem:S-proc}, we direct our attention to proving sufficiency, i.e., the existence of $\tau\colon\mathbb{P}\times\mathbb{P}\rightarrow\mathbb{R}_{>0}$ satisfying~\eqref{eq:VtauW}. To this end, consider $\pazocal{V}$ and $\pazocal{W}$ as in~\eqref{eq:Vcal} and~\eqref{eq:S-proc}, respectively, and suppose that $q^\top\pazocal{V}(\hat{p}_+,\hat{p})q>0$ for all $q\in\mathbb{R}^{n_q}\setminus\{0\}$ such that $q^\top\pazocal{W}(\hat{p}_+,p)q\geqslant 0$. For each $\hat{p}_+,\hat{p}\in\mathbb{P}$, we distinguish two cases: \begin{enumerate*}[label=(\roman*)]
			\item $(H-K(\hat{p}_+,\hat{p})C)\neq 0$ or $K(\hat{p}_+,\hat{p})D\neq 0$, or
			\item $(H-K(\hat{p}_+,\hat{p})C)=0$ and $K(\hat{p}_+,\hat{p})D=0$.		
		\end{enumerate*}
		
		First, we show that $\pazocal{W}(\hat{p}_+,\hat{p})$ satisfies Slater's condition for all $\hat{p}_+,\hat{p}\in\mathbb{P}$ that fall under case (i). To this end, let $\xi=(\bar{e},\frac{1}{2}\Lambda((H-K(\hat{p}_+,\hat{p})C)\bar{e}+K(\hat{p}_+,\hat{p})D\bar{w}),0,\bar{w},0)$, for some $\bar{e}\in\mathbb{R}^{n_x}$ and $\bar{w}\in\mathbb{R}^{n_w}$, then, $2\xi^\top\pazocal{W}(\hat{p}_+,\hat{p})\xi=\|\Lambda((H-K(\hat{p}_+,\hat{p})C)\bar{e}+K(\hat{p}_+,\hat{p})D\bar{w})\|^2$. It follows that, for all $\hat{p}_+,\hat{p}\in\mathbb{P}$ for which $H-K(\hat{p}_+,\hat{p})C\neq 0$ and/or $K(\hat{p}_+,\hat{p})D\neq 0$ (i.e., case (i)), there exist $\bar{e}$ and $\bar{w}$ such that $\xi^\top\pazocal{W}(\hat{p}_+,\hat{p})\xi>0$. Hence, by pointwise application, for each such $\hat{p}_+,\hat{p}\in\mathbb{P}$, of Lemma~\ref{lem:S-proc}, there exists $\bar{\tau}\colon\mathbb{P}\times\mathbb{P}\rightarrow\mathbb{R}_{\geqslant 0}$ such that
		\begin{equation}
			\pazocal{V}(\hat{p}_+,\hat{p})-\bar{\tau}(\hat{p}_+,\hat{p})\pazocal{W}(\hat{p}_+,\hat{p})\succ 0,
		\end{equation}
		for all $\hat{p}_+,\hat{p}\in\mathbb{P}$ for which $H-K(\hat{p}_+,\hat{p})C\neq 0$ and/or $K(\hat{p}_+,\hat{p})D\neq 0$. We now proceed with case (ii), i.e., the pathological case where $H-K(\hat{p}_+,\hat{p})C=0$ and $K(\hat{p}_+,\hat{p})D=0$ for some $\hat{p}_+,\hat{p}\in\mathbb{P}$. In this case Slater's condition clearly does not hold, however, we can use Lemma~\ref{lem:finsler} to show that there exists $\bar{\mu}\colon\mathbb{P}\times\mathbb{P}\rightarrow\mathbb{R}$ such that $\pazocal{V}(\hat{p}_+,\hat{p})-\bar{\mu}(\hat{p}_+,\hat{p})\pazocal{W}(\hat{p}_+,\hat{p})\succ 0$. To see this, note that by definition of $\tilde{\phi}$ (above~\eqref{eq:error-system}) and by~\eqref{eq:prop-slope-restr}, $H-K(\hat{p}_+,\hat{p})C=0$ and $K(\hat{p}_+,\hat{p})D=0$ together imply that $\tilde{\phi}=0$ and, hence $\pazocal{M}^\top q=0$ with $\pazocal{M}$ as defined in~\eqref{eq:M}. By Lemma~\ref{lem:finsler}, there exists $\bar{\mu}\colon\mathbb{P}\times\mathbb{P}\rightarrow\mathbb{R}$ such that $\pazocal{V}(\hat{p}_+,\hat{p})-2\pazocal{M}\pazocal{M}^\top \succ 0$ for all $\hat{p}_+,\hat{p}\in\mathbb{P}$ for which $H-K(\hat{p}_+,\hat{p})C=0$ and $K(\hat{p}_+,\hat{p})D=0$. Since $2\pazocal{M}\pazocal{M}^\top = -\pazocal{W}(\hat{p}_+,\hat{p})$ (when $H-K(\hat{p}_+,\hat{p})C=0$ and $K(\hat{p}_+,\hat{p})D=0$) and since case (i) and (ii) together cover all $\hat{p}_+,\hat{p}\in\mathbb{P}$, we conclude that there exists a function $\tau\colon\mathbb{P}\times\mathbb{P}\rightarrow\mathbb{R}$ such that~\eqref{eq:VtauW} holds for all $\hat{p}_+,\hat{p}\in\mathbb{P}$. Finally, one of the diagonal blocks of~\eqref{eq:VtauW} reads $2\tau(\hat{p}_+,\hat{p})I\succ 0$ and, hence, we have $\tau(\hat{p}_+,\hat{p})>0$ for all $\hat{p}_+,\hat{p}\in\mathbb{P}$, thus, there exists $\tau\colon\mathbb{P}\times\mathbb{P}\rightarrow\mathbb{R}_{>0}$ such that~\eqref{eq:VtauW} holds for all $\hat{p}_+,\hat{p}\in\mathbb{P}$.
\end{pf*}

\begin{pf*}{Proof of Corollary~\ref{cor:param-mismatch-iss}.}
	From either Theorem~\ref{thm:observer-lmi} or Theorem~\ref{thm:suffnec}, we already have ISS with respect to $v$, $w$ and $\psi$. Under the additional assumptions in this corollary, we now show ISS with respect to $\tilde{p}$. Firstly,~\eqref{eq:prop-slope-restr} implies that $\phi$ is globally Lipschitz continuous, i.e., there exists $c\in\mathbb{R}_{\geqslant 0}$ such that $\|\phi(y)-\phi(x)\|\leqslant c\|y-x\|$ for all $x,y\in\mathbb{R}^{n_{\phi}}$. By the boundedness theorem, there exists $K_\phi\in\mathbb{R}_{\geqslant 0}$ such that $\|\phi(x)\|\leqslant K_{\phi}$ for all $x\in\mathbb{B}_{K_x}^{n_x}$. Since $\xi_i$, $i\in\pazocal{N}$, are continuous and $\mathbb{P}$ is compact, there exist $\pazocal{K}$-functions $\alpha_i$, $i\in\pazocal{N}$, such that $\|\xi_i(\tilde{p}+p)-\xi_i(p)\|\leqslant \alpha_i(\|\tilde{p}\|)$ for all $p\in\mathbb{P}$ and $\tilde{p}\in\mathbb{D}\coloneqq\{\hat{p}-p\,|\,p,\hat{p}\in\mathbb{P}\}$~\cite[Corollary III.10]{Lazar2013}. By substitution in $\gamma_{\psi}$ and using the boundedness of $x$, $u$ and $v$, we obtain~\eqref{eq:iss-param} with $\gamma_{\tilde{p}}(s)=\gamma_{\psi}(\alpha_{\tilde{p}}(s))$, where $\alpha_{\tilde{p}}(s)=\sum_{i\in\pazocal{N}}\alpha_i(s)(\|E_i\|+\|A_i\|)K_x+\|B_i\|\Delta_u+\|G_i\|K_\phi+\|F_i\|K_v)$, and $\gamma_v,\gamma_w$ as in Theorem~\ref{thm:observer-lmi} (if the system satisfies Assumption~\ref{ass:Ebar} as in Theorem~\ref{thm:suffnec} we substitute $E_i=\bar{E}$, for all $i\in\pazocal{N}$, in $\alpha_{\tilde{p}}$).
\end{pf*}

\bibliographystyle{plain}        
\bibliography{phd-bibtex}	     
\end{document}